\documentclass[10pt]{article}
\usepackage{amsfonts}
\usepackage[leqno]{amsmath}
\usepackage{graphicx}
\usepackage{latexsym}
\usepackage{amsmath,amsfonts,amssymb,amsthm,mathrsfs,euscript,makeidx,color}
\usepackage{enumerate}
\allowdisplaybreaks[1]
\usepackage{multirow}

\def\sqr#1#2{{\vcenter{\vbox{\hrule height.#2pt
              \hbox{\vrule width.#2pt height#1pt \kern#1pt \vrule width.#2pt}
              \hrule height.#2pt}}}}
\def\signed #1{{\unskip\nobreak\hfil\penalty50
              \hskip2em\hbox{}\nobreak\hfil#1
              \parfillskip=0pt \finalhyphendemerits=0 \par}}
\def\endpf{\signed {$\sqr69$}}

\def\3n{\negthinspace \negthinspace \negthinspace }
\def\2n{\negthinspace \negthinspace }
\def\1n{\negthinspace }
\def\bel{\begin{equation}\label}
\def\eel{\end{equation}}

\def\dbE{\mathbb{E}}
\def\dbF{\mathbb{F}}

\def\dbP{\mathbb{P}}

\def\dbR{\mathbb{R}}
\def\dbS{\mathbb{S}}

\def\dbU{\mathbb{U}}

\def\sE{\mathscr{E}}

\def\sL{\mathscr{L}}

\def\sP{\mathscr{P}}

\def\sR{\mathscr{R}}

\def\sU{\mathscr{U}}

\def\={\buildrel \triangle \over =}

\def\ds{\displaystyle}

\def\ns{\noalign{\ss}}
\def\a{\alpha}
\def\b{\beta }
\def\g{\gamma}
\def\d{\delta}

\def\z{\zeta}

\def\l{\lambda}
\def\m{\mu}
\def\n{\nu}
\def\si{\sigma}

\def\f{\varphi}

\def\i{\infty}
\def\Th{\Theta}
\def\L{\Lambda}

\def\F{\Phi}
\def\O{\Omega}

\def\BTh{{\bf\Theta}}
\def\cE{{\cal E}}
\def\cF{{\cal F}}

\def\BX{{\bf X}}

\def\Bw{{\bf w}}

\def\no{\noindent}

\def\ss{\smallskip}
\def\ms{\medskip}

\def\q{\quad}
\def\qq{\qquad}
\def\hb{\hbox}

\def\limsup{\mathop{\overline{\rm lim}}}
\def\liminf{\mathop{\underline{\rm lim}}}

\def\lan{\mathop{\langle}}
\def\ran{\mathop{\rangle}}

\def\h{\widehat}
\def\wt{\widetilde}

\def\cd{\cdot}
\def\cds{\cdots}

\def\ae{\hbox{\rm a.e.{ }}}
\def\as{\hbox{\rm a.s.}}

\def\diag{\hbox{\rm $\,$diag$\,$}}

\def\deq{\mathop{\triangleq}}
\def\les{\leqslant}
\def\ges{\geqslant}

\def\({\Big (}
\def\){\Big )}
\def\[{\Big[}
\def\]{\Big]}
\def\bde{\begin{definition}\label}
\def\ede{\end{definition}}
\def\be{\begin{equation}}
\def\bel{\begin{equation}\label}
\def\ee{\end{equation}}
\def\bt{\begin{theorem}\label}
\def\et{\end{theorem}}
\def\bc{\begin{corollary}\label}
\def\ec{\end{corollary}}
\def\bl{\begin{lemma}\label}
\def\el{\end{lemma}}
\def\bp{\begin{proposition}\label}
\def\ep{\end{proposition}}
\def\bas{\begin{assumption}}
\def\eas{\end{assumption}}
\def\br{\begin{remark}\label}
\def\er{\end{remark}}
\def\ba{\begin{array}}
\def\ea{\end{array}}
\def\ed{\end{document}}

\def\rf{\eqref}
\newcommand{\ad}{&\!\!\!\displaystyle}

\def\square#1{\vbox{\hrule\hbox{\vrule height#1%
     \kern#1\vrule}\hrule}}
\def\rectangle#1#2{\vbox{\hrule\hbox{\vrule height#1%
     \kern#2\vrule}\hrule}}

\font\tenbb=msbm10 \font\sevenbb=msbm7 \font\fivebb=msbm5

\newfam\bbfam
\scriptscriptfont\bbfam=\fivebb \textfont\bbfam=\tenbb
\scriptfont\bbfam=\sevenbb

\newtheorem{theorem}{Theorem}[section]
\newtheorem{corollary}[theorem]{Corollary}

\newtheorem{lemma}[theorem]{Lemma}
\newtheorem{proposition}[theorem]{Proposition}

\theoremstyle{definition}
\newtheorem{definition}[theorem]{Definition}
\newtheorem{remark}{Remark}

\makeatletter
   
   \@addtoreset{equation}{section}
\makeatother

\begin{document}

\title{\bf Optimal Ergodic Control of Linear Stochastic Differential Equations with Quadratic Cost Functionals Having Indefinite Weights\thanks{This work is supported in part by NSF Grant DMS-1812921, the National Natural Science Foundation of China (11971099), the Science and Technology Development Plan Project of Jilin Province (20190103026JH).}}

\author{Hongwei Mei\footnote{ Department of Mathematics, The University of Kansas, Lawrence, KS 66045, USA; email: {\tt hongwei.mei@ku.edu}},~~
Qingmeng Wei\footnote{School of Mathematics and Statistics, Northeast Normal University, Changchun 130024, China; email: {\tt weiqm100@nenu.} {\tt edu.cn}},~~
Jiongmin Yong\footnote{Department of Mathematics, University of
Central Florida, Orlando, FL 32816, USA; email: {\tt
jiongmin.yong@ucf.edu}} }

\maketitle

\no\bf Abstract: \rm An optimal ergodic control problem (EC problem, for short) is investigated for a linear stochastic differential equation with quadratic cost functional. Constant nonhomogeneous terms, not all zero, appear in the state equation, which lead to the asymptotic limit of the state non-zero. Under the stabilizability condition, for any (admissible) closed-loop strategy, an invariant measure is proved to exist, which makes the ergodic cost functional well-defined and the EC problem well-formulated. Sufficient conditions, including those allowing the weighting matrices of cost functional to be indefinite, are introduced for finiteness and solvability for the EC problem. Some comparisons are made between the solvability of EC problem and the closed-loop solvability of stochastic linear quadratic optimal control problem in the infinite horizon. Regularized EC problem is introduced to be used to obtain the optimal value of the EC problem.

\ms

\no\bf Keywords: \rm  Linear-quadratic problem, ergodic control, indefinite quadratic cost, invariant measure, algebraic Riccati equation.

\ms

\no\bf AMS Mathematics Subject Classification. \rm 93E20, 49N10, 60F17.

\section{Introduction}

Let $(\O,\cF,\dbF,\dbP)$ be a complete filtered probability space on which a standard one-dimensional Brownian motion $\{W(t),t\ges0\}$ is defined such that $\dbF=\{\cF_t\}_{t\ges0}$ is the natural filtration of $W(\cd)$ augmented by all the $\dbP$-null sets in $\cF$. We begin with the following $n$-dimensional controlled linear stochastic differential equation
\bel{state-0}\left\{\2n\ba{ll}
\ds dX(t)=\big[AX(t)+Bu(t)+b\big]dt+\big[CX(t)+Du(t)+\si\big]dW(t),\q t\ges0,\\
\ns\ds X(0)=x.\ea\right.\ee
In the above, $A,C\in\dbR^{n\times n}$, $B,D\in\dbR^{n\times m}$, are called the {\it coefficients} of the system, $b,\si\in\dbR^n$ are called the {\it nonhomogeneous terms}, $X(\cd)$ is the {\it state process} valued in $\dbR^n$, and $u(\cd)$ is the {\it control process} valued in $\dbR^m$. We call \rf{state-0} the {\it state equation}. Let
$$\ba{ll}
\ns\ds\sU[0,T]\1n\equiv\1n\Big\{u:[0,\i)\1n\times\1n\O\to\dbR^m\bigm| u(\cd)\hb{ is $\dbF$-progressively measurable, }\dbE\1n\int_0^T\2n|u(t)|^2dt\1n<\1n\infty\Big\},\ea$$
and
$$\ba{ll}
\ns\ds\sU_{loc}[0,\infty)=\bigcap_{T>0}\sU[0,T],\\
\ns\ds\sU[0,\infty)\equiv L^2_\dbF(0,\infty;\dbR^m)=\Big\{u(\cd)\in\sU_{loc}[0,\infty)\bigm| \dbE\int_0^\infty|u(t)|^2dt<\infty\Big\}.\ea$$
Clearly, for any $u(\cd)\in\sU_{loc}[0,\infty)$ and any {\it initial state} $x\in\dbR^n$, state equation \rf{state-0} admits a unique solution $X(\cd)\equiv X(\cd\,;x,u(\cd))$ which is $\dbF$-adapted and continuous, thus integrable on any finite interval $[0,T]$. To measure the performance of the control $u(\cd)$ on the interval $[0,T]$, we introduce the following {\it cost functional}:
\bel{cost-0}J_T(x;u(\cd))=\dbE\int_0^Tg(X(t),u(t))dt,\ee
where
\bel{g}g(x,u)=\lan Qx,x\ran+2\lan Sx,u\ran+\lan Ru,u\ran+2\lan q,x\ran+2\lan\rho,u\ran,\ee
with some suitable matrices $Q,S,R$ and vectors $q,\rho$. Then an optimal control problem on $[0,T]$ can be formulated:

\ms

\bf Problem (LQ$_{[0,T]}$). \rm For any given $x\in\dbR^n$, find a $\bar u(\cd)\in\sU[0,T]$ such that
\bel{inf[0,T]}J_T(x;\bar u(\cd))=\inf_{u(\cd)\in\sU[0,T]}J_T(x;u(\cd)).\ee

When a $\bar u(\cd)\in\sU[0,T]$ exists satisfying \rf{inf[0,T]}, we say that Problem (LQ$_{[0,T]}$) is {\it open-loop solvable} and $\bar u(\cd)$ is called an {\it open-loop optimal control}, the corresponding state process $\bar X(\cd)\equiv X(\cd\,;x,\bar u(\cd))$ is called the corresponding {\it open-loop optimal state process} and $(\bar X(\cd),\bar u(\cd))$ is called an {\it open-loop optimal pair}. Problem (LQ$_{[0,T]}$) is referred to as an LQ problem on $[0,T]$.

\ms

It is well-known by now that under proper conditions, Problem (LQ$_{[0,T]}$) (even allowing $b(\cd),\si(\cd),q(\cd),\rho(\cd)$ to be random) admits a unique open-loop optimal control $\bar u(\cd)\in\sU[0,T]$ which has a {\it closed-loop representation} via the solution to a Riccati differential equation; Further, this coincides with the outcome of a so-called {\it closed-loop optimal strategy} (see, for example, \cite{Sun-Li-Yong 2016}, for details). It is natural to ask what it will be if we consider the problem on $[0,\infty)$? Namely, consider the same state equation \rf{state-0} with the following cost functional:
\bel{cost-1}J_\infty(x;u(\cd))=\dbE\int_0^\infty g(X(t),u(t))dt.\ee
Such LQ problems have been studied in \cite{Sun-Yong 2018} (see the references cited therein as well for some details) with $b,\si,q,\rho$ replaced by globally square integrable $\dbF$-progressively measurable processes $b(\cd),\si(\cd),q(\cd),\rho(\cd)$ on $[0,\infty)$ and the homogeneous system, denoted by $[A,C;B,D]$ is {\it stabilizable}, by which we mean that there exists a matrix $\Th\in\dbR^{m\times n}$, called a {\it stabilizer} of $[A,C;B,D]$, such that the homogeneous {\it closed-loop system}
\bel{closed-loop-0}\left\{\2n\ba{ll}
\ds dX_0(t)= (A+B\Th)X_0(t) dt+ (C+D\Th)X_0(t) dW(t),\qq t\ges0,\\
\ns\ds X_0(0)=x\ea\right.\ee
admits a unique solution $X_0(\cd)\in L_\dbF^2(0,\infty;\dbR^n)$. Now, in the case that $b,\si$ are constant vectors, not all are zero, the (global) integrability condition is not satisfied. Consequently, even $[A,C;B,D]$ is stabilizable, the state $X(t;x,u(\cd))$ will not approach to zero as $t\to\i$. Thus, the corresponding cost functional will not be well-defined in general. Hence, the corresponding LQ problem is not well-formulated in the traditional way.

\ms

In this paper, we are going to formulate an LQ problem with the state equation \rf{state-0} and a quadratic cost functional which is closely related to the original \rf{cost-1}. We will develop a theory for that kind of LQ problems. In order our procedure can be carried out, throughout the paper, we will assume the following.

\ms

{\bf(H1)} The matrices $A,C\in\dbR^{n\times n}$, $B,D\in\dbR^{n\times m}$ satisfy the following:
\bel{K ne}\BTh[A,C;B,D]=\big\{\Th\in\dbR^{m\times n}\bigm|\Th\hb{ stabilizes }[A,C;B,D]\big\}\ne\varnothing.\ee

When (H1) holds, we call any pair $(\Th,v)\in\BTh[A,C;B,D]\times\dbR^m$ an {\it admissible closed-loop strategy} (see \cite{Sun-Yong 2018} for a similar notion), and define
\bel{U}\dbU=\big\{u:\dbR^n\to\dbR^m\bigm|u(x)=\Th x+v,\q (\Th,v)\in\BTh[A,C;B,D]\times\dbR^m\big\},\ee
which is the set of all {\it outcomes} of the admissible closed-loop strategies. Any $u(\cd)\in\dbU$ is also called a {\it linear feedback control}, or a {\it closed-loop control}. In what follows, we will identify $u(x)\equiv\Th x+v$ with $(\Th,v)$, via which, one has
\bel{U=Th}\dbU=\BTh[A,C;B,D]\times\dbR^m.\ee
It is clear that for any $(\Th,v)\in\dbU$, the following {\it closed-loop system}:
\bel{closed-loop-1}\left\{\1n\ba{ll}
\ds dX(t)=\big[(A+B\Th)X(t)+Bv+b\big]dt+\big[(C+D\Th)X(t)+Dv+\si\big]dW(t),\q t\ges0,\\
\ns\ds X(0)=x,\ea\right.\ee
has a unique solution $X(\cd)\equiv X(\cd\,;x,\Th,v)$ on $[0,\infty)$. Although it is not necessarily in $L^2_\dbF(0,\infty;\dbR^n)$, we will show (in the next section) that the following holds:
\bel{sup|X|}\sup_{t\ges0}\dbE|X(t;x,\Th,v)|^2<\infty.\ee
Hence, for any $\l>0$, the so-called {\it Abel mean} type functional can be defined:
\bel{cost-2}J^\l_\infty(x;\Th,v)=\dbE\int_0^\infty e^{-\l t}g\big(X(t;x,\Th,v),\Th X(t;x,\Th,v)+v\big)dt,\qq(\Th,v)\in\dbU.\ee
Consequently, one could try to find a $(\bar\Th_\l,\bar v_\l)\in\dbU$ such that
\bel{inf Je}J_\infty^\l(x;\bar\Th_\l,\bar v_\l)=\inf_{(\Th,v)\in\dbU}J_\infty^\l(x;\Th,v).\ee
It is natural to further ask what happen if we send $\l\to0^+$ (trying to recover the solution to the original problem in some sense)? Since under (H1), only \rf{sup|X|} is guaranteed, we could not expect the limit $\ds\lim_{\l\to0^+}J_\infty^\l(x;\Th,v)$ to exist and to be finite. It turns out that the following exists:
\bel{E}\wt J_\infty(x;\Th,v)=\liminf_{\l\to 0^+}\l J^\l_\infty(x;\Th,v)\deq\cE(\Th,v),\qq\forall(\Th,v)\in\dbU,\ee
with some function $\cE:\dbU\to\dbR$, independent of $x\in\dbR^n$, called an {\it ergodic cost function}. Hence, we could formulate the following optimal control problem.

\ms

\bf Problem (EC). \rm Find a pair $(\bar\Th,\bar v)\in\dbU$ such that
\bel{min cE}\cE(\bar\Th,\bar v)=\inf_{(\Th,v)\in\dbU}\cE(\Th,v).\ee

\ms

We call the above an {\it ergodic optimal control problem}. Any pair $(\bar\Th,\bar v)\in\dbU$ satisfying the above is called an {\it optimal strategy} of Problem (EC). Note that when this happens, we have
\bel{exp}J_\i^\l(x;\bar\Th,\bar v)={1\over\l}\cE(\bar\Th,\bar v)+o(\frac1\l),\qq\l\to0^+.\ee

\ms

Note that, in general, $\cE(\cd\,,\cd)$ could be complicated. Therefore, directly approach to such a problem is difficult. Hence, we would like to find an equivalent form which could be easier to handle. To this end, we take a different viewpoint. For any $(\Th,v)\in\dbU$ and $x\in\dbR^n$, the solution $X(\cd)\equiv X(\cd\,;x,\Th,v)$ of \rf{closed-loop-1} is a homogeneous Markov process. By \rf{sup|X|}, $\{X(t)\bigm|t\ges0\}$ is tight. Moreover, $X(\cd\,;x,\Th,v)$ is a {\it Feller process}, by which we mean that $x\mapsto\dbE\big[h(X(t;x,\Th,v)\big]$ is continuous for any bounded continuous function $h(\cd)$ and $t\ges0$. Hence, by \cite{Khas 2011} and taking into account the linearity of the state equation, we see that $X(\cd\,;x,\Th,v)$ admits a unique {\it invariant measure} $\pi^{\Th,v}$, indicating the dependence on $(\Th,v)\in\dbU$. That is to say if the initial state $X(0)$ follows the distribution $\pi^{\Th,v}$, then for each $t>0$, $X(t)$ follows the same distribution $\pi^{\Th,v}$. In next section, we will show that
\bel{abel-er}\wt J_\infty(x;\Th,v)\equiv\liminf_{\l\to 0^+}\l J^\l_\infty(x;\Th,v)\equiv\cE(\Th,v)=\int_{\dbR^n}g(x,\Th x+v)\pi^{\Th,v}(dx),\qq\forall x\in\dbR^n.\ee

In this paper, we are going to investigate Problem (EC). The main novelty of this paper can be briefly summarized as follows:

\ms

(i) Under the condition that the homogeneous system $[A,C;B,D]$ is stabilizable, we rigorously formulate the ergodic optimal control problem by means of invariant measure.

\ms

(ii) The finiteness and solvability of Problem (EC) will be discussed and sufficient conditions for these notions will be obtained, in terms of algebraic Riccati equation/inequality.

\ms

(iii) Comparison will be made between the results of Problem (EC) and classical LQ problem in the infinite horizon $[0,\i)$. It turns out that the algebraic Riccati equation for the solvability of Problem (EC) is the same as that for the closed-loop solvability of the classical LQ problem. Hence, to ensure the solvability of Problem (EC), one could just look at the closed-loop solvability of Problem (LQ), which is further equivalent to the open-loop solvability of LQ problem. However we point out that the solvability of the classical LQ problem in an infinite horizon is not necessary for that of Problem (EC).

\ms

(iv) We find a general sufficient condition (not just present some particular examples) for the uniform convexity of the cost functional for the {\it stabilized} LQ problem. The condition allows either $Q$ or $R$ to be negative to some extent. This combining the above (iii), we obtain a set of sufficient conditions for the solvability of Problem (EC).

\ms

(v) When Problem (EC) is merely finite, by introducing the regularized EC problem, we find a general scheme to find the optimal value of Problem (EC).

\ms

Study of deterministic LQ problems can be traced back to the works of Bellman--Glicksberg--Gross \cite{Bellman-Glicksberg-Gross 1958}, Kalman \cite{Kalman 1960}, and Letov \cite{Letov 1961} in the later 1950s and early 1960s. Investigation of stochastic LQ problems was initiated by Wonham \cite{Wonham 1968} in 1968. See \cite{Bismut 1976, Anderson-Moore 1989} and references cited therein for some other followed-up works. In all these classical works, the following classical {\it positive definiteness condition}
\bel{QR>0}R>0,\q Q-S^\top R^{-1}S\ges0\ee
has been taken granted for a long time. In 1977, Molinari found that for deterministic LQ problem, $Q$ could be a little negative (see also \cite{You 1983} for a more general case in Hilbert spaces). In 1998, Chen--Li--Zhou \cite{Chen-Li-Zhou 1998} further found that for stochastic LQ problem, even $R$ could be a little negative, see \cite{Yong-Zhou 1999, Ait Rami-Zhou 2000, Chen-Yong 2001, Ait Rami-Chen-Moore-Zhou 2001, Ait Rami-Moore-Zhou 2002,
Hu-Zhou 2003, Li-Zhou-Ait Rami 2003, Huang-Li-Yong 2015,Sun-Li-Yong 2016, Sun-Yong 2018}, for many further development.

\ms

On the other hand, the ergodic control problem for general stochastic diffusion rather than linear-quadratic ones, has been investigated in the book \cite{Arapostathis-Borkar-Kumar 2014} (see also the references cited therein). The main approach  is to analyze  the so-called stationary Hamilton-Jacobi equation (see Section 3.6.2 there). To guarantee the well-posedness of the stationary HJB equations, most of the results in \cite{Arapostathis-Borkar-Kumar 2014} requires that the diffusion of the system is non-degenerate and the cost functional is inf-compact (or called near-monotone in
\cite{Arapostathis-Borkar-Kumar 2014}).
%For linear-quadratic cases, the stationary Hamilton-Jacobi equation is equivalent to an algebraic Riccati equation. Thus we might investigate a larger class of ergodic control problems for linear-quadratic systems beyond definite assumptions of cost functions.
For the case $Q>0$ and $R=I$, the linear-quadratic ergodic control problem with stationary (random) coefficients was studied in \cite{Guatteri-Masiero 2009,Guatteri-Masiero 2009b}.

\ms

The rest of the paper is arranged as follows. In Section 2, we present some preliminary results, including the derivation of the ergodic cost function by means of invariant measure. Section 3 is devoted to the study of Problem (EC). Sufficient conditions will be obtained for the finiteness and the solvability of the problem, respectively. In Section 4, we will make a formal comparison between the solvability of Problem (EC) and the closed-loop solvability of the corresponding classical LQ problem in the infinite horizon. Also, a sufficient condition is introduced so that the cost functional of the stabilized LQ problem is uniformly convex with respect to the control. This will lead to the solvability of Problem (EC). Section 5 is concerned with the finiteness of Problem (EC). The optimal value of the cost function is obtained via the regularized ergodic problem. General one-dimensional situation is considered in Section 6. Final concluding remarks are collected in Section 7.

\section{Preliminary}

In this section, we will present some preliminary results. We introduce some spaces.
$$L^p_\cF(\O;\dbR^n)=\big\{\xi:\O\to\dbR^n\bigm|\xi\hb{ is $\cF$-measurable, }\dbE|\xi|^p<\infty\big\},\qq p\ges1.$$
$$C_b(\dbR^n)=\big\{h:\dbR^n\to\dbR\bigm|h(\cd)\hb{ is continuous and uniformly bounded}\big\}.$$

\subsection{Estimates of the state process}

In this subsection, we will briefly discuss the stabilization of the homogeneous system  $[A,C;B,D]$, and its consequences. For any $\Th\in\BTh[A,C;B,D]$, the homogeneous closed-loop system \rf{closed-loop-0} admits a unique solution $X_0(\cd)\in L^2_\dbF(0,\infty;\dbR^n)$. Denote
$$F(\Th)=(A+B\Th)+(A+B\Th)^\top+(C+D\Th)^\top (C+D\Th)\in\dbS^n,$$
where $\dbS^n$ is the set of all $(n\times n)$ symmetric (real) matrices. Then by It\^o's formula, we obtain
$${d\over dt}\(\dbE|X_0(t)|^2\)=\dbE\lan F(\Th)X_0(t),X_0(t)\ran,\q t\ges0.$$
There exists an orthogonal matrix $\F\equiv\F(\Th)$ such that
$$\F^\top F(\Th)\F=\L\equiv\diag(\l_1,\cds,\l_n),\q\l_1\ges\l_2\ges\cds\ges\l_n,$$
with $\l_1,\cds,\l_n$ being all the eigenvalues of $F(\Th)$. Consequently, by letting $\wt X_0=\F^\top X_0$, we have
$${d\over dt}\dbE|\wt X_0(t)|^2=\dbE\lan\L\wt X_0(t),\wt X_0(t)\ran.$$
Hence, by $X_0(\cd)\in L^2_\dbF(0,\infty;\dbR^n)$, we have $\wt X_0(\cd)\in L^2_\dbF(0,\infty;\dbR^n)$, and it is necessary that $\l_k<0$, $k=1,2,\cds,n$. If we denote
$$-\l(\Th)=\max\{\l_1,\cds,\l_n\}\equiv\max\si\big(F(\Th)\big),$$
then
\bel{}\dbE|X_0(t)|^2\les K_\Th e^{-\l(\Th)t}|x|^2,\qq\forall t\ges0.\ee
Here, the constant $K_\Th>0$ depends on $\Th$ through $\F(\Th)$, and $\l(\Th)>0$ depends on $\Th$ and intrinsically depends on $[A,C;B,D]$. We also note from the above that
\bel{}\lan F(\Th)x,x\ran=\lan\F\L\F^\top x,x\ran\les-\l(\Th)|\F^\top x|^2=-\l(\Th)|x|^2.\ee
The following lemma is concerned with the estimate \rf{sup|X|} and beyond.

\bl{Lemma-2.1} \sl For any strategy $(\Th,v)\in\dbU$ and $x\in\dbR^n$, the closed-loop system \rf{closed-loop-1} admits a unique solution $X(\cd)=X(\cd\,;x,\Th,v)$ such that
\bel{sup|X|*}\sup_{t\ges0}\dbE|X(t)|^2\les K(1+|x|^2).\ee
Hereafter, $K>0$ stands for a generic constant which can be different from line to line. Moreover, let $(\h\Th,\h v)\in\dbU$, $\h x\in\dbR^n$, and let $\h X(\cd)$ be the corresponding solution of \rf{closed-loop-1}, then
\bel{|BX|}\dbE|X(t)-\h X(t)|^2\les e^{-{[\l(\Th)\vee\l(\h\Th)]\over2}t}|x-\h x|^2+K\(|v-\h v|^2+|\Th-\h\Th|^2\),\qq\forall t\ges0,\ee
for some $K=K_{\Th,v,\h\Th,\h v}>0$, which is continuous in $(\Th,v,\h\Th,\h v)$.

\el

\it Proof. \rm Let $(\Th,v)\in\dbU$ and $x\in\dbR^n$. It is standard that the closed-loop system \rf{closed-loop-1} admits a unique solution $X(\cd)\equiv X(\cd\,;x,\Th,v)$. By  It\^o's formula, we have
\bel{L2bounded}\ba{ll}
\ns\ds{d\over dt}\(\dbE|X(t)|^2\)=\dbE\(\lan F(\Th)X(t),X(t)\ran\\
\ns\ds\qq\qq\qq\qq+2\lan(Bv+b)+\sum_{k=1}^d(C_k+D_k\Th)^\top (D_kv+\si_k),X(t)\ran+\sum_{k=1}^d|D_kv+\si_k|^2\)\\
\ns\ds\les\dbE\(-\l(\Th)|X(t)|^2+2\lan(Bv+b)+\sum_{k=1}^d(C_k+D_k\Th)^\top (D_kv+\si_k),X(t)\ran\)+\sum_{k=1}^d|D_kv+\si_k|^2\\
\ns\ds\les-{\l(\Th)\over2}\dbE|X(t)|^2\2n+{2\over\l(\Th)}\Big|Bv\1n+\1n b+\1n\sum_{k=1}^d(C_k\1n
+\1n D_k\Th)^\top (D_kv\1n+\1n\si_k)\Big|^2\2n+\1n\sum_{k=1}^d|D_kv\1n+\1n\si_k|^2\\
\ns\ds=-{\l(\Th)\over2}\dbE|X(t)|^2+L_0(\Th),\ea\ee
with
$$L_0(\Th)={2\over\l(\Th)}\Big|Bv+b+\sum_{k=1}^d(C_k+D_k\Th)^\top (D_kv+\si_k)\Big|^2+\sum_{k=1}^d|D_kv+\si_k|^2.$$
Hence,
\bel{|X|}\dbE|X(t)|^2\les e^{-\l(\Th)t\over2}|x|^2+\int_0^te^{-{\l(\Th)(t-s)\over2}}L_0(\Th)ds=e^{-\l(\Th)t\over2}|x|^2
+{2L_0(\Th)\over\l(\Th)}\(1-e^{-\l(\Th)t\over2}\).\ee
Consequently, we obtain from the above that
$$\sup_{t\in[0,\infty)}\dbE|X(t)|^2\les\(|x|^2\vee
{2L_0(\Th)\over\l(\Th)}\).$$
Thus, \rf{sup|X|*} follows.

\ms

Next, let $\BX(\cd)=X(\cd)-\h X(\cd)$. Then the following holds:
\bel{closed-loop-2}\left\{\1n\ba{ll}
\ds d\BX(t)=\big[(A+B\Th)\BX(t)+B(v-\h v)+B(\Th-\h\Th)\h X(t)\big]dt\\
\ns\ds\qq\qq+\sum_{k=1}^d\big[(C_k+D_k\Th)\BX(t)+D_k(v-\h v)+D_k(\Th-\h\Th)\h X(t)\big]dW_k(t),\q t\ges0,\\
\ns\ds\BX(0)=x-\h x.\ea\right.\ee
Similar to the above, we see that
$${d\over dt}\(\dbE|\BX(t)|^2\)\les-{\l(\Th)\over2}\dbE|\BX(t)|^2+L_0(t;\Th,\h \Th),$$
with
$$\ba{ll}
\ns\ds L_0(t;\Th,\h\Th)={2\over\l(\Th)}\Big|B(v-\h v)+B(\Th-\h\Th)\h X(t)
+\sum_{k=1}^d(C_k+D_k\Th)^\top\big[D_k(v-\h v)+D_k(\Th-\h\Th)\h X(t)\big]\Big|^2\\
\ns\ds\qq\qq\qq\qq+\sum_{k=1}^d|D_k(v-\h v)+D_k(\Th-\h\Th)\h X(t)|^2\\
\ns\ds\qq\les{K\over\l(\Th)}(1+|\Th|^2)\(|v-\h v|^2+|\h X(t)|^2|\Th-\h\Th|^2\)\\
\ns\ds\qq\les{K\over\l(\Th)}(1+|\Th|^2)\[|v-\h v|^2+\(|\h x|^2\vee{2L_0(\h\Th)\over\l(\h\Th)}\)|\Th-\h\Th|^2\]\equiv L_0(\Th,\h\Th).\ea$$
By Gronwall's inequality,
$$\ba{ll}
\ns\ds\dbE|X(t)-\h X(t)|^2\les e^{-\l(\Th)t\over2}|x-\h x|^2+\int_0^te^{-{\l(\Th)(t-s)\over2}}L_0(\Th,\h\Th)ds\\
\ns\ds=e^{-\l(\Th)t\over2}|x-\h x|^2
+{2L_0(\Th,\h\Th)\over\l(\Th)}\(1-e^{-\l(\Th)t\over2}\).\ea$$
By exchange the position of $X(\cd)$ and $\h X(\cd)$, we finally obtain \rf{|BX|}. \endpf

\ms

%From the above proof, we see that the set $\BTh[A,C;B,D]$ is open.

\subsection{Invariant Measures}

In this subsection, we will review some results on invariant measures. For any Euclidean space $\dbR^\ell$, let $\sL(\dbR^\ell)$ be its Lebesgue $\si$-field, and
$$\sP_2(\dbR^\ell)=\Big\{\n:\sL(\dbR^\ell)\to[0,1]\bigm|\n\hb{ is a probability on $\sL(\dbR^\ell)$, }
\int_{\dbR^\ell}|x|^2\n(dx)<\infty\Big\}.$$
For $\m_1,\m_2\in\sP_2(\dbR^n)$, we define
\bel{w_2}\ba{ll}
\ns\ds \Bw_2(\m_1,\m_2)=\inf\Big\{\(\int_{\dbR^{2n}}|x_1-x_2|^2\n(dx_1,dx_2)\)^{1\over2}\bigm|
\n\in\sP_2(\dbR^{2n}),\\
\ns\ds\qq\qq\qq\qq\qq\qq\n(dx_1,\dbR^n)=\m_1(dx_1),~\n(\dbR^n,dx_2)=\m_2(dx_2)\Big\}.\ea\ee
The above is called the {\it Wasserstein-2 metric} (or simply {\it $\Bw_2$-metric}), under which $\sP_2(\dbR^n)$ is a complete metric space (see Theorem 6.16 in \cite{Villani 2008}). For a random variable $\xi$, we denote $\text{law}(\xi)$ to be the distribution of $\xi$. By the definition of $\Bw_2$, we have
\bel{w2L2}\Bw_2^2(\text{law}(\xi),\text{law}(\eta))\les\dbE|\xi-\eta|^2,\qq\forall\xi,\eta\in L^2_\cF(\O;\dbR^n).\ee
The following proposition gives an equivalent condition of convergence under Wasserstein-2 metric (see \cite{Villani 2008}).

\begin{proposition}\label{equiw2} \sl Let $\m_k,\m\in\sP_2(\dbR^n)$. Then
$$\lim_{k\to\infty}\Bw_2(\m_k,\m)=0,$$
if and only if $\m_k$ weakly converges to $\m$, i.e.,
$$\lim_{k\to\infty}\int_{\dbR^n}h(x)\m_k(dx)=\int_{\dbR^n}h(x)\m(dx),\qq\forall h(\cd)\in C_b(\dbR^n),$$
and also
$$\lim_{k\to\infty}\int_{\dbR^n}|x|^2\m_k(dx)=\int_{\dbR^n}|x|^2\m(dx).$$
\end{proposition}
\ms

We know that for any closed-loop strategy $(\Th,v)\in\dbU$, the unique solution map $x\mapsto X(\cd\,;x)$ of \rf{closed-loop-1} is a stochastic flow (\cite{Kunita 1990}) which can be uniquely characterized by its {\it transition probability} $p(t,x;dy)$, where
$$p(t,x;dy)=\dbP\big(X(t;x)\in dy\big),\qq(t,x)\in[0,\i)\times\dbR^n.$$
%
%Actually, $p(t,x;\cd)$ is nothing but the law of the random variable $X(t;x)$.
We have the following lemma.

\bl{inv} \sl For any closed-loop strategy $(\Th,v)\in\dbU$, let $p(t,x;\cd)$ be the transition probability of the stochastic flow $X(t;x)$ of \eqref{closed-loop-1}. Then there exists a unique invariant measure $\pi$ such that
\bel{lim}\lim_{t\to\infty}\Bw_2\big(p(t,x;\cd),\pi\big)=0,\qq\forall x\in\dbR^n.\ee
Moreover, if $(\Theta_k,v_k)\in\BTh[A,C;B,D]\times\dbR^n$ converges to some $(\Th,v)\in\BTh[A,C;B,D]\times\dbR^n$, then $\pi^{\Th_k,v_k}$ converges to $\pi^{\Th,v}$ in $\Bw_2$-metric.
	
\end{lemma}

\begin{proof} \rm We want to show that given any $x$, $\{p(t,x;\cd):t\ges0\}$ is Cauchy, as $t\to\infty$ in $(\sP_2(\dbR^n),\Bw_2)$ with a same limit for any $x$. To prove this, we let $\Psi$ be the set of couples $(\f,\psi)$ such that $\f$ and $\psi$ are bounded continuous with $\f(y_1)+\psi(y_2)\les|y_1-y_2|^2$. Using the Kantorovich's duality (see Theorem 5.9 in \cite{Villani 2008}), for $t_2>t_1\ges0$, and $x_1,x_2\in\dbR^n$, we have
$$\ba{ll}
\ns\ad\Bw^2_2\big(p(t_1,x_1;\cd), p(t_2,x_2;\cd)\big)\\
\ns\ad=\sup_{(\f,\psi)\in\Psi}\(\int_{\dbR^n}\f(y_1)p(t_1,x_1;dy_1)+\int_{\dbR^n}\psi(y_2)
p(t_2,x_2;dy_2)\)\\
\ns\ad=\sup_{(\f,\psi)\in\Psi}\(\int_{\dbR^n}\f(y_1)p(t_1,x_1;dy_1)+\int_{\dbR^n}\psi(y_2)
\int_{\dbR^n}p(t_1,z;dy_2)p(t_2-t_1,x_2;dz)\)\\
\ns\ad\les\int_{\dbR^n}p(t_2-t_1,x_2;dz)\[\sup_{(\f,\psi)\in\Psi}\(\int_{\dbR^n}\f(y_1)
p(t_1,x_1;dy_1)+\int_{\dbR^n}\psi(y_2)p(t_1,z;dy_2)\)\]\\
\ns\ad\les\int_{\dbR^n}p(t_2-t_1,x_2;dz)\Bw^2_2\big(p(t_1,x_1;\cd),p(t_1,z;\cd)\big)\\
\ns\ad\les\int_{\dbR^n}\dbE|X(t_1;x_1)-X(t_1;z)|^2p(t_2-t_1,x_2;dz)\\
\ns\ad\les e^{-{\l(\Th)\over2}t_1}\int_{\dbR^n}|x_1-z|^2p(t_2-t_1,x_2;dz).\ea$$
In the last two steps, we have used \rf{w2L2} and \eqref{|BX|}. By letting $t_2>t_1\to\infty$, we see that $\{p(t,x;\cd):t\ges0\}$ is Cauchy (as $t\to\infty$) in $(\sP_2(\dbR^n),\Bw_2)$ with some limit $\pi$. Note that $\pi$ is an invariant measure which is independent of the choice of $x$. Moreover, if $\pi'$ is another invariant measure, then for any $h\in C_b(\dbR^n)$, one has
$$\int_{\dbR^n}h(y)\pi'(dy)=\int_{\dbR^n}h(y)\int_{\dbR^n}p(t,z;dy)\pi'(dz)\to
\int_{\dbR^n}h(y)\pi(dy)\pi'(dz)=\int_{\dbR^n}h(y)\pi(dy).$$
In the second equality, we take $t\rightarrow \infty$. This proves that $\pi=\pi'$. Thus the  invariant measure $\pi$ is unique.

\ms

Finally, if $(\Theta_k,v_k)\in\BTh[A,C;B,D]\times\dbR^n$ converges to some $(\Th,v)\in\BTh[A,C;B,D]\times\dbR^n$, by \eqref{|BX|},  we can see that $\Bw_2(\pi^{\Th_k,v_k},\pi^{\Th,v})\rightarrow 0$. The proof is complete.
\end{proof}

\subsection{The ergodic cost functional}

In this subsection, we will prove the claim \eqref{abel-er}.

\ms

Let $(\Th,v)\in\dbU$ be fixed, and let $X(\cd)\equiv X(\cd\,;\xi,\Th,v)$ be the solution of the closed-loop system \rf{closed-loop-1} with the initial state $\xi\in\dbR^n$. We introduce the following:
\bel{nu}\n^\xi_\l(G)\equiv\l\dbE\int_0^\infty e^{-\l s}I(X(s;\xi)\in G)ds,\qq G\in\cF,\ee
which is called the {\it occupation measure} of $X(\cd\,;\xi)$.
Then, with $u(x)=\Th x+v$, we have
\bel{}\l J^\l_\infty(\xi;u(\cd))=\l\int_0^\infty e^{-\l t}\dbE g\big(X(t;\xi),u(X(t;\xi))\big)dt=\int_{\dbR^n}g(x,u(x))\n^\xi_\l(dx).\ee
We want to prove that $\n^\xi_\l(dx)$ converges to $\pi^u$ weakly as $\l\to0$.

\ms

Note that under closed-loop strategy $(\Th,v)$, $X(\cd\,;\xi)$ is a homogeneous Markov process. Thus, we may let $q^\xi(t,x;dy)$ be its transition probability, i.e.,
$$q^\xi(t,x;dy)=\dbP\big(X(s+t;\xi)\in dy,X(s;\xi)=x\big).$$
Now, for any continuous bounded function $f:\dbR^n\to\dbR$, one has
$$\ba{ll}
\ns\ds\int_{\dbR^n}f(y)\int_{\dbR^n}q^\xi(t,x;dy)\n^\xi_\l(dx)\\
\ns\ds=\int_{\dbR^n}f(y)\int_{\dbR^n}q^\xi(t,x;dy)\l\dbE\int_0^\infty e^{-\l s}I(X(s;\xi)\in dx)ds\\
\ns\ds=\l\int_{\dbR^n}f(y)\dbE\int_0^\infty e^{-\l s}\int_{\dbR^n}q^\xi(t,x;dy)I(X(s;\xi)\in dx)ds\\
\ns\ds=\l\int_{\dbR^n}f(y)\dbE\int_0^\infty e^{-\l s}I(X(t+s;\xi)\in dy)ds\\
\ns\ds=\l\int_{\dbR^n}f(y)e^{\l t}\dbE\int_0^\infty e^{-\l s}I(X(s;\xi)\in dy)ds-\l e^{\l t}\int_{\dbR^n}f(y)\dbE\int_0^t e^{-\l s}I(X(s;\xi)\in dy)ds\\
\ns\ds=e^{\l t}\int_{\dbR^n}f(y)\n_\l^\xi(dy)-\l e^{\l t}\int_{\dbR^n}f(y)\dbE\int_0^te^{-\l s}I(X(s;\xi)\in dy)ds.\ea$$
For any fixed $t>0$, letting $\l\to0^+$, we see that the second term on the right-hand side will go to zero. Since $\n_\l^\xi$ is tight (because $X(t;\xi)$ is tight), any subsequence has a weakly convergent subsequence with a same limit $\pi$. Note that $x\mapsto\int_{\dbR^n}f(y)q^\xi(t,x;dy)$ is continuous, by the Feller property, then for any $t>0$,
$$\int_{\dbR^n}f(y)\int_{\dbR^n}q^\xi(t,x;dy)\pi(dx)=\int_{\dbR^n}f(y)\pi(dy).$$
This verifies that $\pi$ is an invariant measure. By the uniqueness of the invariant measure $\pi^u$, $\pi=\pi^u$. This shows that $\n^\xi_\l$ converges to $\pi^u$ weakly. Note that by Proposition \ref{equiw2} and Lemma \ref{inv}, it follows
$$\int_{\dbR^n}|x|^2\n_\l^\xi(dx)=\l\int_0^\infty e^{-\l s}\dbE|X(t;\xi)|^2ds\to\int_{\dbR^n}|x|^2\pi^u(dx),\q\text{as $\l\to0^+$}.$$
As a result,
$$\wt J_\infty(x;u(\cdot))=\lim_{\l\to0^+}\l J^\l_\infty(\xi;u(\cd))=\lim_{\l\to0^+}\int_{\dbR^n}g(x,u(x))\n^\xi_\l(dx)
=\int_{\dbR^n}g(x,u(x))\pi^u(dx).$$
This verifies our claim \eqref{abel-er}.

\ms

Under (H1), for any $u(\cd)\in\dbU$, we may also introduce the following so-called {\it Ces\`aro mean} type cost functional:
\bel{Cesaro}\wt J_T(x;u(\cd))={1\over T}J_T(x;u(\cd))\equiv{1\over T}\int_0^Tg(X(t),u(t))dt.\ee
If we introduce the following corresponding occupation measure
$$\wt\n^u_T(dx)={1\over T}\dbE\int_0^TI(X(t)\in dx)dt,$$
then, with a similar argument (details can be found in Theorem 3.1.1 of \cite{Da Prato-Zabczyk 1996}), one has
\bel{Cesaro2}\liminf_{T\to\infty}{1\over T}J_T(x;u(\cd))=\liminf_{T\to\infty}\int_{\dbR^n}g(x,u(x))\n_T^u(dx)=\int_{\dbR^n}g(x,u(x))
\pi^u(dx).\ee
%
%In this paper, we will use Abel-mean cost in our proof instead of the average-type cost \rf{Cesaro2}.

\section{Ergodic Optimal Control Problem --- Finiteness and Solvability}

In this section, we investigate the  ergodic optimal control problem. For convenience, let us recall the problem as follows.

\ss

\bf Problem (EC). \rm Let (H1) hold. Find a $\bar u(\cd)\in\dbU$ such that
\bel{inf*-1}\cE(\bar u(\cd))=\inf_{u(\cd)\in\dbU}\cE(u(\cd))\equiv\sE.\ee

\ms

Note that under (H1), $\dbU\ne\varnothing$. Hence, there will be at least one strategy $u(\cd)\in\dbU$ such that $\cE(u(\cd))$ is finite, which implies $\sE<\infty$. Adopting the usual terminology of optimal LQ problems, we introduce the following definition.

\bde{} \rm Problem (EC) is said to be {\it finite} if $\sE>-\infty$. If there (uniquely) exists a $\bar u(\cd)\in\dbU$ satisfying \rf{inf*-1}, Problem (EC) is said to be (uniquely) {\it solvable}. In this case, $\bar u(\cd)$ is called an (the) {\it optimal strategy} of
Problem (EC).

\ede
For simplicity, if $u(x)=\Theta x+v$, we also write  $\cE(\Theta,v)\equiv\cE(u(\cdot))$.
Recall that
\bel{cE}\cE(\Th,v)=\int_{\dbR^n}g(x,\Th x+v)\pi^{\Th,v}(dx),\qq\forall(\Th,v)\in\BTh[A,C;B,D]
\times\dbR^n.\ee
Therefore, in the case that
\bel{g>-K}g(x,u)\equiv\Bigg\langle\begin{pmatrix}Q&S^\top\\ S&R\end{pmatrix}\begin{pmatrix}x\\ u\end{pmatrix},\begin{pmatrix}x\\ u\end{pmatrix}\Bigg\rangle+2\Bigg\langle\begin{pmatrix}q\\ \rho\end{pmatrix}\begin{pmatrix}x\\ u\end{pmatrix}\Bigg\rangle\ges-K,\qq\forall(x,u)\in\dbR^n\times\dbR^m,\ee
for some $K\ges0$, one will have
\bel{E>-K}\cE(\Th,v)\ges-K,\qq\forall(\Th,v)\in\BTh[A,C;B,D]\times\dbR^n,\ee
leading to the finiteness of Problem (EC). Note that \rf{g>-K} is equivalent to the following:
\bel{Q>0}\begin{pmatrix}Q&S^\top\\ S&R\end{pmatrix}\ges0,\qq\begin{pmatrix}q\\ \rho\end{pmatrix}\in\sR\Bigg(\begin{pmatrix}Q&S^\top\\ S&R\end{pmatrix}\Bigg).\ee
We refer to the above as the classical {\it positive semi-definiteness condition}. Apparently, condition \rf{Q>0} is too restrictive. As a matter of fact, by assuming \rf{Q>0}, one does not make use of the compatibility of $g(x,\Th x+v)$ and the related invariant measure $\pi^{\Th,v}(\cd)$. On the other hand, we recall that in standard stochastic LQ theory (\cite{Chen-Li-Zhou 1998, Sun-Li-Yong 2016, Sun-Xiong-Yong 2020}), $Q$ or $R$ is even allowed to be a little negative (therefore \rf{Q>0} fails) within a certain extent, still keeping the corresponding problem to have optimal controls. This inspires us to explore the possible relaxation on \rf{Q>0} below.

\ms

Note that for any $(\Th,v)\in\BTh[A,C;B,D]\times\dbR^n$, one has
\bel{g*}\ba{ll}
\ns\ds g(x,\Th x+v)=\Bigg\langle\begin{pmatrix}Q&S^\top\\ S&R\end{pmatrix}\begin{pmatrix}x\\ \Th x+v\end{pmatrix},\begin{pmatrix}x\\ \Th x+v\end{pmatrix}\Bigg\rangle+2\Bigg\langle\begin{pmatrix}q\\ \rho\end{pmatrix}\begin{pmatrix}x\\ \Th x+v\end{pmatrix}\Bigg\rangle\\
\ns\ds=\lan Qx,x\ran+2\lan Sx,\Th x+v\ran+\lan R(\Th x+v),\Th x+v\ran+2\lan q,x\ran+2\lan\rho,\Th x+v\ran\\
\ns\ds=\lan(Q+S^\top\Th+\Th^\top S+\Th^\top R\Th)x,x\ran+2\lan(S+R\Th)x,v\ran+\lan Rv,v\ran+2\lan q+\Th^\top\rho,x\ran+2\lan\rho,v\ran\\
\ns\ds=\Bigg\langle\begin{pmatrix}S^\top\Th+\Th^\top S+Q&(S+R\Th)^\top\\ S+R\Th&R\end{pmatrix}\begin{pmatrix}x\\ v\end{pmatrix},\begin{pmatrix}x\\ v\end{pmatrix}\Bigg\rangle+2\Bigg\langle\begin{pmatrix}q+\Th^\top\rho\\ \rho\end{pmatrix}\begin{pmatrix}x\\ v\end{pmatrix}\Bigg\rangle.\ea\ee
Thus,
\bel{cE}\ba{ll}
\ns\ds\cE(\Th,v)=\int_{\dbR^n}\[\lan(Q+S^\top\Th+\Th^\top S+\Th^\top R\Th)x,x\ran+2\lan(S+R\Th)x,v\ran+\lan Rv,v\ran\\
\ns\ds\qq\qq\qq\qq+2\lan q+\Th^\top\rho,x\ran+2\lan\rho,v\ran\]\pi^{\Th,v}(dx).\ea\ee
Now, we would like to find another representation of function $\cE(\cd)$, which will help us to obtain the finiteness and solvability of Problem (EC). To this end, let us make some preparations.

\ms

For any $\Pi\in\dbS^{n\times n}$ and $\Th\in\BTh[A,C;B,D]$, we denote
$$L_\Pi=B^\top\Pi+D^\top\Pi C+S,\qq M_{\Th,\Pi}=\begin{pmatrix}
Q_{\Th,\Pi}&L_\Pi^\top+\Th^\top(R+D^\top\Pi D)\\
L_\Pi+(R+D^\top\Pi D)\Th&R+D^\top\Pi D\end{pmatrix},$$
with
\bel{Lyapunov}\ba{ll}
\ns\ds Q_{\Th,\Pi}=\Pi(A+B\Th)+(A+B\Th)^\top\Pi+(C+D\Th)^\top\Pi(C+D\Th)+S^\top\Th+\Th^\top S+\Th^\top R\Th+Q\\
\ns\ds\qq~=\Pi A+A^\top \Pi+C^\top\Pi C+L_\Pi^\top\Th+\Th^\top L_\Pi+\Th^\top(R+D^\top\Pi D)\Th+Q.\ea\ee
If $\Pi\in\dbS^n$ such that
\bel{R>0*}R+D^\top\Pi D\ges0,\qq\sR(B^\top\Pi+D^\top\Pi+S)\subseteq\sR(R+D^\top\Pi D),\ee
then there exists a $\L_\Pi\in\dbR^{m\times n}$ such that $L_\Pi=(R+D^\top\Pi D)\L_\Pi$ which leads to the following:
$$\ba{ll}
\ns\ds L_\Pi^\top\Th+\Th^\top L_\Pi+\Th^\top(R+D^\top\Pi D)\Th\\
\ns\ds=\L_\Pi^\top(R+D^\top\Pi D)\Th+\Th^\top(R+D^\top\Pi D)\L_\Pi+\Th^\top(R+D^\top\Pi D)\Th\\
\ns\ds=(\L_\Pi+\Th)^\top(R+D^\top\Pi D)(\L_\Pi+\Th)-\L_\Pi^\top(R+D^\top\Pi D)(R+D^\top\Pi D)^\dag(R+D^\top\Pi D)\L_\Pi\\
\ns\ds=\big[\Th+(R+D^\top\Pi D)^\dag L_\Pi\big]^\top(R+D^\top\Pi D)\big[\Th+(R+D^\top\Pi D)^\dag L_\Pi\big]-L_\Pi^\top(R+D^\top\Pi D)^\dag L_\Pi.\ea$$
Hence,
\bel{}\ba{ll}
\ns\ds Q_{\Th,\Pi}=\Pi A+A^\top \Pi+C^\top\Pi C+Q-L_\Pi^\top(R+D^\top\Pi D)^\dag L_\Pi+(\L_\Pi+\Th)^\top(R+D^\top\Pi D)(\L_\Pi+\Th)\\
\ns\ds\qq~\equiv\h Q_\Pi+\big[\Th+(R+D^\top\Pi D)^\dag L_\Pi\big]^\top(R+D^\top\Pi D)\big[\Th+(R+D^\top\Pi D)^\dag L_\Pi\big],\ea\ee
with
\bel{h Q}\ba{ll}
\ns\ds\h Q_\Pi=\Pi A+A^\top \Pi+C^\top\Pi C+Q-L_\Pi^\top(R+D^\top\Pi D)^\dag L_\Pi\\
\ns\ds\qq\equiv\Pi A+A^\top \Pi+C^\top\Pi C+Q-(B^\top\Pi+D^\top\Pi C+S)^\top(R+D^\top\Pi D)^\dag(B^\top\Pi+D^\top\Pi C+S).\ea\ee
Also, when \rf{R>0*} holds, one has
\bel{M>}\ba{ll}
\ns\ds M_{\Th,\Pi}=\begin{pmatrix}I&(\L_\Pi\1n+\1n\Th)^\top\\0&I\end{pmatrix}\begin{pmatrix}
Q_{\Th,\Pi}\1n-\1n(\L_\Pi\1n+\1n\Th)^\top(R\1n+\1nD^\top\1n\Pi D)(\L_\Pi\1n+\1n\Th)&0\\0&R\1n+\1n D^\top\Pi D\end{pmatrix}\begin{pmatrix}I&0\\ \L_\Pi+\Th&I\end{pmatrix}\\
\ns\ds\qq\q=\begin{pmatrix}I&(\L_\Pi\1n+\1n\Th)^\top\\0&I\end{pmatrix}\begin{pmatrix}
\h Q_\Pi&0\\0&R\1n+\1n D^\top\Pi D\end{pmatrix}\begin{pmatrix}I&0\\ \L_\Pi\1n+\1n\Th&I\end{pmatrix}.\ea\ee
Consequently, in the case that \rf{R>0*} holds and the following {\it algebraic Riccati inequality} holds
\bel{Q>0*}\h Q_\Pi\equiv\Pi A+A^\top \Pi+C^\top\Pi C+Q-(B^\top\Pi+D^\top\Pi C+S)^\top(R+D^\top\Pi D)^\dag(B^\top\Pi+D^\top\Pi C+S)\ges0,\ee
one has $M_{\Th,\Pi}\ges0$. Further, if we let
\bel{Th_0}\ba{ll}
\ns\ds\Th_0=-(R+D^\top\Pi D)^\dag(B^\top\Pi+D^\top\Pi C+S)+[I-(R+D^\top\Pi D)^\dag(R+D^\top\Pi D)]\L\\
\ns\ds\q~\equiv-(R+D^\top\Pi D)^\dag L_\Pi+[I-(R+D^\top\Pi D)^\dag(R+D^\top\Pi D)]\L,\ea\ee
for any $\L\in\dbR^{m\times n}$, then, noting $L_\Pi=(R+D^\top\Pi D)\L_\Pi$, we have
\bel{RTh}(R+D^\top\Pi D)\Th_0=-(R+D^\top\Pi D)(R+D^\top\Pi D)^\dag(R+D^\top\Pi D)\L_\Pi=-L_\Pi.\ee
Hence,
\bel{Th>Th}\ba{ll}
\ns\ds Q_{\Th,\Pi}-Q_{\Th_0,\Pi}=L_\Pi^\top(\Th-\Th_0)+(\Th-\Th_0)^\top L_\Pi+\Th^\top(R+D^\top\Pi D)\Th-\Th_0^\top(R+D^\top\Pi D)\Th_0\\
\ns\ds\qq~=L_\Pi^\top(\Th-\Th_0)+(\Th-\Th_0)^\top L_\Pi+(\Th-\Th_0)^\top(R+D^\top\Pi D)\Th_0+\Th_0^\top(R+D^\top\Pi D)(\Th-\Th_0)\\
\ns\ds\qq\qq\qq+(\Th-\Th_0)^\top(R+D^\top\Pi D)(\Th-\Th_0)\\
\ns\ds\qq~=(\Th-\Th_0)^\top(R+D^\top\Pi D)(\Th-\Th_0)\ges0,\qq\forall\Th\in\dbR^{m\times n}.\ea\ee
In another word, $\Th_0$ defined by \rf{Th_0} is a minimum of the map $\Th\mapsto Q_{\Pi,\Th}$, taking the usual order in $\dbS^n$. For given $\Pi\in\dbS^n$, the set of all $\Th_0$ of form \rf{Th_0} is denoted by $\Upsilon[\Pi]$, i.e.,
\bel{Upsilon}\Upsilon[\Pi]=\Big\{-(R+D^\top\Pi D)^\dag(B^\top\Pi+D^\top\Pi C+S)+[I-(R+D^\top\Pi D)^\dag(R+D^\top\Pi D)]\L\bigm|\L\in\dbR^{m\times n}\Big\}.\ee
In the case that $R+D^\top\Pi D>0$, $\Upsilon[\Pi]$ is a singleton. Also, we see that for any $\Th_0\in\Upsilon[\Pi]$, noting \rf{RTh},
\bel{}\ba{ll}
\ns\ds Q_{\Th_0,\Pi}=\Pi A+A^\top\Pi+C^\top\Pi C+Q+L_\Pi^\top\Th_0+\Th_0^\top L_\Pi+\Th_0^\top(R+D^\top\Pi D)\Th_0\\
\ns\ds\qq\q=\Pi A+A^\top\Pi+C^\top\Pi C+Q-\Th_0(R+D^\top\Pi D)\Th_0\\
\ns\ds\qq\q=\Pi A+A^\top \Pi+C^\top\Pi C+Q-L_\Pi^\top(R+D^\top\Pi D)^\dag L_\Pi=\h Q_\Pi.\ea\ee
Next, we note that for any $\Th\in\BTh[A,C;B,D]$, we know that system $[A+B\Th,C+D\Th]$ is asymptotically stable. Therefore,
$$A+B\Th+(A+B\Th)^\top+(C+D\Th)^\top(C+D\Th)<0,$$
which leads to the invertibility of $A+B\Th$. We now ready to present the following result.

\bl{revalue} \sl Let {\rm(H1)} hold. For any $(\Th,v)\in\BTh[A,C;B,D]\times\dbR^n$, let $\pi^{\Th,v}$ be the corresponding invariant measure. Then for any $\Pi\in\dbS^n$, the ergodic cost function $\cE(\cd)$ admits the following representation:
\bel{rep1}\cE(\Th,v)\1n=\2n\int_{\dbR^n}\2n\Bigg\langle M_{\Th,\Pi}\begin{pmatrix}x\\v\end{pmatrix},\begin{pmatrix}x\\v\end{pmatrix}\Bigg\rangle\pi^{\Th,v}(dx)+2\lan B^\top \eta_{_{\Th,\Pi}}+D^\top \Pi\si+\rho,v\ran+\lan\Pi\si,\si\ran+2\lan\eta_{_{\Th,\Pi}},b\ran,\ee
where $\eta_{_{\Th,\Pi}}\in\dbR^n$ is the solution to the following linear equation:
\bel{eta}(A+B\Th)^\top\eta_{_{\Th,\Pi}}+\Pi b+(C+D\Th)^\top\Pi\si+q+\Th^\top\rho=0.\ee
\el

\begin{proof} For $\Th\in\BTh[A,C;B,D]$, let $\pi^{\Th,v}$ be the invariant measure. Then if we let $X(\cd\,;\xi)$ be the solution of \rf{closed-loop-1} with the initial state $\xi$ having the distribution $\pi^{\Th,v}(\cd)$. Then for any $t>0$, $X(t)$ will have the same distribution $\pi^{\Th,v}(\cd)$. Thus, for any $(\Pi,\eta)\in\dbS^n\times\dbR^n$, $t\mapsto\dbE\big[\lan\Pi X(t),X(t\ran+2\lan\eta,X(t)\ran\big]$ stays as a constant. Hence, by It\^o's formula,  we obtain ($t$ will be suppressed)
$$\ba{ll}
\ns\ds0={d\over dt}\(\dbE\lan\Pi X(t),X(t)\ran+2\dbE\lan\eta,X(t)\ran\)\\
\ns\ns\q=\dbE\lan\Pi[(A+B\Th)X+Bv+b],X\ran+\dbE\lan\Pi X,(A+B\Th)X+Bv+b\ran\\
\ns\ds\qq+\dbE\lan\Pi[(C+D\Th)X+Dv+\si],(C+D\Th)X+Dv+\si\ran+2\dbE\lan\eta, (A+B\Theta)X+Bv+b\ran\\
\ns\ds\q=\dbE\lan\big[\Pi(A+B\Th)+(A+B\Th)^\top\Pi+(C+D\Th)^\top\Pi(C+D\Th)\big]X,X\ran\\
\ns\ds\qq+2\dbE\lan\big[B^\top\Pi+D^\top\Pi(C+D\Th)\big]X,v\ran+2\dbE\lan\Pi b+(C+D\Th)^\top\Pi\si+(A+B\Th)^\top\eta,X\ran\\
\ns\ds\qq+\lan D^\top\Pi Dv,v\ran+2\dbE\lan B^\top\eta+D^\top \Pi\si,v\ran+\lan\Pi\si,\si\ran+2\lan\eta,b\ran\\
\ns\ds\q=\dbE\lan\big[\Pi A+A^\top\Pi+C^\top\Pi C+(\Pi B+C^\top\Pi D)\Th+\Th^\top(B^\top\Pi
+D^\top\Pi C)+\Th^\top D^\top\Pi D\Th\big]X,X\ran\\
\ns\ds\qq+2\dbE\lan\big[B^\top\Pi+D^\top\Pi(C+D\Th)\big]X,v\ran+2\dbE\lan\Pi b+(C+D\Th)^\top\Pi\si+(A+B\Th)^\top\eta,X\ran\\
\ns\ds\qq+\lan D^\top\Pi Dv,v\ran+2\dbE\lan B^\top\eta+D^\top \Pi\si,v\ran+\lan\Pi\si,\si\ran+2\lan\eta,b\ran\\
\ns\ds\q=\dbE\lan Q_{\Th,\Pi}\1n-\1n(S^\top\1n\Th\1n+\1n\Th^\top\1n S\1n+\1n
\Th^\top\1n R\Th\1n+\1n Q)X,X\ran
-\lan Rv,v\ran+2\dbE\lan\big[B^\top\Pi+D^\top\Pi(C+D\Th)\big]X,v\ran\\
\ns\ds\qq+\lan(R+D^\top\Pi D)v,v\ran+2\dbE\lan\Pi b+(C+D\Th)^\top\Pi\si+(A+B\Th)^\top \eta,X\ran +2\lan D^\top \Pi\si+B^\top\eta,v\ran\\
\ns\ds\qq+\lan\Pi\si,\si\ran+2\lan\eta,b\ran\\
\ns\ds\q=\dbE\lan Q_{\Th,\Pi}X,X\ran-\cE(\Th,v)+2\dbE\lan\big[L_\Pi+(R+D^\top\Pi D)\Th\big]X,v\ran\\
\ns\ds\qq+2\dbE\lan\Pi b+(C+D\Th)^\top\Pi\si+(A+B\Th)^\top \eta+q+\Th^\top\rho,X\ran\\
\ns\ds\qq+\lan (R+ D^\top\Pi D)v,v\ran+2\lan D^\top \Pi\si+B^\top \eta+\rho,v\ran+\lan\Pi\si,\si\ran+2\lan\eta,b\ran.\ea$$
This implies that
\bel{cE=*}\ba{ll}
\ns\ds\cE(\Th,v)=\dbE\lan Q_{\Th,\Pi}X,X\ran+2\dbE\big\langle\big[L_\Pi+(R+D^\top\Pi D)\Th\big]X,v\big\rangle\\
\ns\ds\qq\qq\q+2\dbE\lan\Pi b+(C+D\Th)^\top\Pi\si+(A+B\Th)^\top \eta+q+\Th^\top\rho,X\ran\\
\ns\ds\qq\qq\q+\lan (R+ D^\top\Pi D)v,v\ran+2\lan D^\top \Pi\si+B^\top \eta+\rho,v\ran+\lan\Pi\si,\si\ran+2\lan\eta,b\ran.\ea\ee
Taking $\eta=\eta_{_{\Th,\Pi}}$, we have
$$\ba{ll}
\ns\ds\cE(\Th,v)=\dbE\lan Q_{\Th,\Pi}X,X\ran+2\dbE\big\langle\big[L_\Pi
+(R+D^\top\Pi D)\Th\big]X,v\ran+\lan(R+ D^\top\Pi D)v,v\big\langle\\
\ns\ds\qq\qq\qq+2\lan B^\top \eta_{_{\Th,\Pi}}+D^\top \Pi\si+\rho,v\ran+\lan\Pi\si,\si\ran+2\lan\eta_{_{\Th,\Pi}},b\ran\\
\ns\ds\qq\qq=\int_{\dbR^n}\[\lan Q_{\Th,\Pi}x,x\ran+2\big\langle\big[L_\Pi
+(R+D^\top\Pi D)\Th\big]x,v\big\rangle+\lan(R+ D^\top\Pi D)v,v\ran\]\pi^{\Th,v}(dx)\\
\ns\ds\qq\qq\qq+2\lan B^\top \eta_{_{\Th,\Pi}}+D^\top \Pi\si+\rho,v\ran+\lan\Pi\si,\si\ran+2\lan\eta_{_{\Th,\Pi}},b\ran.\ea$$
This completes the proof. \end{proof}

\ms

Next, we present a finiteness and solvability theorem for Problem (EC), recalling \rf{Upsilon} for the definition of $\Upsilon[\Pi_0]$.

\bt{finite-solvable} \sl Let {\rm(H1)} hold.

\ms

{\rm(i)} Let $\Pi_0\in\dbS^n$ solve the following algebraic Riccati inequality
\bel{ARE-0}\left\{\2n\ba{ll}
\ds\Pi_0A\1n+\1n A^\top\1n\Pi_0\1n+\1n C^\top\1n\Pi_0C\1n+\1n Q\1n-\1n(B^\top\1n\Pi_0
\1n+\1n D^\top\1n\Pi_0C\1n+\1n S)^\top\1n(R\1n+\1n D^\top\Pi_0D)^\dag\1n(B^\top\Pi_0\1n+\1n D^\top\1n\Pi_0C\1n+\1nS)\1n\ges\1n0,\\
\ns\ds R+D^\top\Pi_0D\ges0,\qq\sR(B^\top\Pi_0+D^\top\Pi_0C+S)\subseteq\sR(R+D^\top\Pi_0D).\ea\right.\ee
Let $(\Th_0,\eta_0)\in\Upsilon[\Pi_0]\times\dbR^n$ such that
\bel{H3equi}\left\{\2n\ba{ll}
\ds B^\top\eta_0+D^\top\Pi_0\si+\rho\in\sR(R+D^\top\Pi_0D),\\
\ns\ds(A+B\Th_0)^\top\eta_0+\Pi_0b+(C+D\Th_0)^\top\Pi_0\si+q+\Th^\top_0\rho\in\sR(Q_{\Th_0, \Pi_0}).\ea\right.\ee
Then $\cE(\cd)$ is bounded from below uniformly on $\dbU$ and Problem (EC) is finite.

\ms

{\rm(ii)} Let $\Pi_0\in\dbS^n$ solve the following algebraic Riccati equation:
\bel{ARE-1}\left\{\2n\ba{ll}
\ds\Pi_0A\1n+\1n A^\top\1n\Pi_0\1n+\1n C^\top\1n\Pi_0C\1n+\1n Q\1n-\1n(B^\top\1n\Pi_0
\1n+\1n D^\top\1n\Pi_0C\1n+\1n S)^\top\1n(R\1n+\1n D^\top\Pi_0D)^\dag\1n(B^\top\Pi_0\1n+\1n D^\top\1n\Pi_0C\1n+\1nS)\1n=\1n0,\\
\ns\ds R+D^\top\Pi_0D\ges0,\qq\sR(B^\top\Pi_0+D^\top\Pi_0C+S)\subseteq\sR(R+D^\top\Pi_0D).\ea\right.\ee
Let $(\bar\Th_0,\bar\eta_0)\in\big\{\Upsilon[\Pi_0]\cap\BTh[A,C;B,D]\big\}\times\dbR^n$ such that
\bel{H3equi-2}\left\{\2n\ba{ll}
\ds B^\top\bar\eta_0+D^\top\Pi_0\si+\rho\in\sR(R+D^\top\Pi_0 D),\\
\ns\ds(A+B\bar\Th_0)^\top\bar\eta_0+\Pi_0b+(C+D\bar\Th_0)^\top\Pi_0\si+q+\bar\Th^\top_0\rho=0.
\ea\right.\ee
then Problem (EC) is solvable with $(\bar\Th_0,\bar v_0)$ being an optimal strategy, where
\bel{v_0}\bar v_0=-(R+D^\top\Pi_0 D)^\dag(D^\top\Pi_0\si+B^\top\bar\eta_0+\rho)+\big[I-(R+D^\top\Pi_0 D)^\dag(R+D^\top\Pi_0 D)^\dag\big]\n,\ee
for any $\n\in\dbR^m$.
\et

\begin{proof} (i) Taking $\Pi=\Pi_0$, $\eta=\eta_0$ in \rf{cE=*},
making use of \rf{Th>Th}, and noting $\Th_0\in\Upsilon[\Pi_0]$ (having property
\rf{RTh}),
\bel{bounbelowca}\ba{ll}
\ns\ds\cE(\Th,v)=\dbE\lan Q_{\Th,\Pi_0}X,X\ran+2\dbE\lan\big[B^\top\Pi_0+D^\top\Pi_0C+S+(R+D^\top\Pi_0 D)\Th\big]X,v\ran\\
\ns\ds\q + 2\dbE\lan\Pi_0 b+(C+D\Theta)^\top\Pi_0\si+(A+B\Th)^\top \eta_0+q+\Th^\top\rho,X\ran\\
\ns\ds\q+\lan (R+ D^\top\Pi_0D)v,v\ran+2\lan D^\top\Pi_0\si+B^\top \eta_0+\rho,v\ran+\lan\Pi_0\si,\si\ran+2\lan\eta_0,b\ran\\
\ns\ds=\dbE\lan Q_{\Th_0,\Pi_0}X,X\ran+\dbE\lan\big[(\Th-\Th_0)^\top(R+D^\top\Pi_0 D)(\Th-\Th_0)\big]X,X\ran\\
\ns\ds\q+2\dbE\lan(D^\top \Pi_0 D+R)(\Th-\Th_0)X,v\ran\\
\ns\ds\q+2\dbE\lan(A+B\Th)^\top \eta_0+\Pi_0 b+(C+D\Th)^\top\Pi_0\si+q+\Th^\top\rho,X\ran\\
\ns\ds\q+\lan (R+D^\top\Pi_0D)v,v\ran+2\lan D^\top\Pi_0\si+B^\top \eta_0+\rho,v\ran+\lan\Pi_0\si,\si\ran+2\lan\eta_0,b\ran\\
\ns\ds=\dbE\lan Q_{\Th_0,\Pi_0}X,X\ran+2\dbE\lan(A+B\Th_0)^\top\eta_0+\Pi_0 b+(C+D\Th_0)^\top\Pi_0\si+q+\Th^\top_0\rho,X\ran\\
\ns\ds\q+\dbE\lan(R\1n+\1n D^\top\1n\Pi_0 D)\big[(\Th\1n-\1n\Th_0)X\1n+\1n v\big],(\Th\1n-\1n\Th_0)X\1n+\1n v\ran\1n
+2\dbE\lan D^\top\Pi_0\si\1n+\1n B^\top\1n\eta_0\1n+\1n\rho,(\Th\1n-\1n\Th_0)X\1n+\1n v\ran\\
\ns\ds\q+\lan\Pi_0\si,\si\ran+2\lan\eta_0,b\ran.\ea\ee
Now, by our assumption, $Q_{\Th_0,\Pi_0}\ges0$, together with \rf{H3equi}, we see that
\bel{cE**}\ba{ll}
\ns\ds\cE(\Th,v)=\dbE\Big|Q_{\Th_0,\Pi_0}^{1\over2}X+[Q_{\Th_0,\Pi_0}^\dag]^{1\over2}
\big\{(A+B\Th_0)^\top\eta_0+\Pi_0 b+(C+D\Th_0)^\top\Pi_0\si+q+\Th^\top_0\rho\big\}\Big|^2\\
\ns\ds\qq\qq\q+\dbE\Big|(R+D^\top\Pi_0D)^{1\over2}\big[(\Th-\Th_0)X+v\big]
+[(R+D^\top\Pi_0D)^\dag]^{1\over2}\big\{D^\top\Pi_0\si\1n+\1n B^\top\1n\eta_0\1n+\1n\rho\big\}\Big|^2\\
\ns\ds\qq\qq\q-\Big|\big[Q_{\Th_0,\Pi_0}^\dag\big]^{1\over2}\big\{(A+B\Th_0)^\top\eta_0+\Pi_0 b+(C+D\Th_0)^\top\Pi_0\si+q+\Th^\top_0\rho\big\}\Big|^2\\
\ns\ds\qq\qq\q-\Big|\big[(R+D^\top\Pi_0D)^\dag\big]^{1\over2}\big\{D^\top\Pi_0\si\1n+\1n B^\top\1n\eta_0\1n+\1n\rho\big\}\Big|^2+\lan\Pi_0\si,\si\ran+2\lan\eta_0,b\ran.\ea\ee
By dropping the first two terms on the right-hand side, we obtain (i).

\ms

(ii) We point out that in the above, the choice of $(\Pi_0,\eta_0)$ does not change the value of $\cE(\Th,v)$. Now, for the current case, we take $(\Th_0,\eta_0)=(\bar\Th_0,\bar\eta_0)$ in \rf{cE**} with $\bar\eta_0$ being the solution to the second equation in \rf{H3equi-2} and note $Q_{\bar\Th_0,\Pi_0}=0$. Then \rf{cE**} becomes
\bel{cE***}\ba{ll}
\ns\ds\cE(\Th,v)=\dbE\Big|(R+D^\top\Pi_0D)^{1\over2}\big[(\Th-\bar\Th_0)X+v
+(R+D^\top\Pi_0D)^\dag(D^\top\Pi_0\si\1n+\1n B^\top\1n\bar\eta_0\1n+\1n\rho)\big]\Big|^2\\
\ns\ds\qq\qq\q-\Big|\big[(R+D^\top\Pi_0D)^\dag\big]^{1\over2}\big\{D^\top\Pi_0\si\1n+\1n B^\top\1n\bar\eta_0\1n+\1n\rho\big\}\Big|^2+\lan\Pi_0\si,\si\ran+2\lan\bar\eta_0,b\ran\\
\ns\ds\qq\qq\ges-\Big|\big[(R+D^\top\Pi_0D)^\dag\big]^{1\over2}\big\{D^\top\Pi_0\si\1n+\1n B^\top\1n\bar\eta_0\1n+\1n\rho\big\}\Big|^2+\lan\Pi_0\si,\si\ran+2\lan\bar\eta_0,b\ran
=\cE(\bar\Th_0,\bar v_0),\ea\ee
proving (ii). \end{proof}

\ms

For later convenience, we introduce the following.

\ms

{\bf(H2)} Let (H1) and \rf{ARE-0}--\rf{H3equi} hold for some $\Pi_0\in\dbS^n$ and some $(\Th_0,\eta_0)\in\Upsilon[\Pi_0]\times\dbR^n$.

\ms

{\bf(H3)} Let (H1) and \rf{ARE-1}--\rf{H3equi-2} hold for some $\Pi_0\in\dbS^n$ and some $(\bar\Th_0,\bar\eta_0)\in\big\{\Upsilon[\Pi_0]\cap\BTh[A,C;B,D]\big\}\times\dbR^n$.

\ms

Then, according to Theorem \ref{finite-solvable}, we have that Problem (EC) is finite if (H2) holds and solvable if (H3) holds. The following corollary is concerned with the classical positive-definite case.

\bc{} \sl Let {\rm(H1)} hold and
\bel{QR>0}\begin{pmatrix}Q&S^\top\\S&R\end{pmatrix}>0.\ee
Then Problem (EC) admits an optimal strategy $(\h\Th,\h v)$ given by the following:
$$\left\{\2n\ba{ll}
\ds\h\Th=-(R+D^\top\h\Pi D)^{-1}(B^\top\h\Pi+D^\top\h\Pi C+S)\in\BTh[A,C;B,D],\\
\ns\ds\h\eta=-\big[(A+B\h\Th)^\top\big]^{-1}[\h\Pi b+(C+D\h\Th)^\top\h\Pi\si+q+\h\Th^\top\rho],\\
\ns\ds\h v=-(R+D^\top\h\Pi D)^{-1}(B^\top\h\eta+D^\top\h\Pi\si+\rho),\ea\right.$$
where $\h\Pi$ is the solution to the following algebraic Riccati equation:
$$\h\Pi A+A^\top\h\Pi+C^\top\h\Pi C+Q-(B^\top\h\Pi+D^\top\h\Pi C+S)^\top(R+D^\top\h\Pi D)^{-1}(B^\top\h\Pi+D^\top\h\Pi C+S)=0.$$
In this case,
$$\sE=\cE(\h\Th,\h v)=\lan\h\Pi\si,\si\ran+2\lan\h\eta,b\ran-\lan(R+D^\top\h\Pi D)\h v,\h v\ran.$$

\ec	

\section{Comparison with Standard LQ Problems}

In this section, we recall some results on classical LQ problems in the infinite horizon $[0,\i)$, with certain improvements, and make some comparisons between these results and results of ergodic problems presented in the previous section.

\ms

Consider the following linear controlled SDE:
\bel{LQ-state}\left\{\2n\ba{ll}
\ds dX(t)=[AX(t)+Bu(t)+b(t)]dt+[CX(t)+Du(t)+\si(t)]dW(t),\qq t\ges0,\\
\ns\ds X(0)=x,\ea\right.\ee
with the cost functional
\bel{LQ-cost}\ba{ll}
\ns\ds\h J_\infty(x;u(\cd))=\dbE\int_0^\infty\(\lan QX(t),X(t)\ran+2\lan SX(t),u(t)\ran+\lan Ru(t),u(t)\ran\\
\ns\ds\qq\qq\qq\qq\qq+2\lan q(t),X(t)\ran+2\lan\rho(t),u(t)\ran\)dt.\ea\ee
Let (H1) hold and assume in addition that
\bel{L^2}b(\cd),\si(\cd),q(\cd)\in L^2_\dbF(0,\i;\dbR^n),\q\rho(\cd)\in L^2_\dbF(0,\i;\dbR^m).\ee
Then the following set of admissible controls is non-empty:
\bel{U_ad}\sU_{ad}[0,\i)=\big\{u(\cd)\in\sU[0,\i)\bigm|X(\cd\,;x,u(\cd))\in L^2_\dbF(0,\i;\dbR^n)\big\},\ee
and $\h J_\infty(x;u(\cd))$ is well-defined for each $u(\cd)\in\sU_{ad}[0,\i)$. Therefore, the following ({\it nonhomogeneous}) LQ problem on $[0,\i)$ is well-formulated.

\ms

{\bf Problem (LQ)$_\i$}. For given $x\in\dbR^n$, find a $\bar u(\cd)\in\sU_{ad}[0,\i)$
such that
\bel{open-loop}\h J_\infty(x;\bar u(\cd))=\inf_{u(\cd)\in\sU_{ad}[0,\i)}\h J_\infty(x;u(\cd)).\ee

\ms

Any $\bar u(\cd)\in\sU_{ad}[0,\infty)$ satisfying \rf{open-loop} is called an {\it open-loop optimal control}. When such a control exists, we say that Problem (LQ)$_\i$ is {\it open-loop solvable} at $x\in\dbR^n$. If Problem (LQ)$_\i$ is open-loop solvable at every $x\in\dbR^n$, we simply say that the problem is open-loop solvable.

\ms

An important special case is that
\bel{b=0}b(\cd)=\si(\cd)=q(\cd)=0,\q\si(\cd)=0.\ee
When the above holds, the problem is said to be {\it homogeneous}. We denote the corresponding state process by $X^0(\cd)=X^0(\cd\,;x,u(\cd))$, the cost functional by
$\h J^0_\infty(x;u(\cd))$, and the problem by Problem (LQ)$^0_\i$. It is not hard to see that the admissible control sets $\sU_{ad}[0,\i)$ for any nonhomogeneous problems (with condition \rf{L^2}), including the homogeneous one, are the same.

\ms

\bde{} \rm (i) Any element $(\Th,v(\cd))\in\dbU[0,\i)\equiv\BTh[A,C;B,D]\times\sU[0,\i)$ is called a {\it closed-loop strategy} of Problem (LQ)$_\i$.

\ms

(ii) Problem (LQ)$_\i$ is said to be {\it closed-loop} solvable if there exists
a $(\bar\Th,\bar v(\cd))\in\dbU[0,\i)$ such that
\bel{}\h J_\i(x;\bar\Th,\bar v(\cd))=\inf_{(\Th,v(\cd))\in\dbU[0,\i)}\h J_\i(x;\Th,v(\cd)).\ee

\ede

The following result is essentially found in \cite{Sun-Yong 2018}.

\bt{Sun-Yong} \sl Let {\rm(H1)} and \rf{L^2} hold. Then Problem {\rm(LQ)$_\i$} is closed-loop solvable if and only if the following algebraic Riccati equation admits a solution $P\in\dbS^n$:
\bel{ARE1}\left\{\2n\ba{ll}
\ds PA+A^\top P+C^\top PC+Q-(B^\top P+D^\top PC+S)^\top(R+D^\top PD)^\dag(B^\top P+D^\top
PC+S)=0,\\
\ns\ds R+D^\top PD\ges0,\qq\sR(B^\top P+D^\top PC+S)\subseteq\sR(R+D^\top PD),\ea\right.\ee
such that for some $\L\in\dbR^{m\times n}$,
\bel{in Th}\ba{ll}
\ns\ds-(R+D^\top PD)^\dag(B^\top P+D^\top PC+S)+\big[I-(R+D^\top PD)^\dag(R+D^\top PD)\big]\L\in\BTh[A,C;B,D],\ea\ee
and the following backward stochastic differential equation (BSDE) admits an adapted solution $(\eta(\cd),\z(\cd))\in L^2_\dbF(0,\i;\dbR^n)^2$:
\bel{BSDE1}\ba{ll}
\ns\ds d\eta(t)=-\Big\{\big[A-B(R+D^\top PD)^\dag(B^\top P+D^\top PC+S)\big]^\top\eta(t)\\
\ns\ds\qq\qq\q+\big[C-D(R+D^\top PD)^\dag(B^\top P+D^\top PC+S)\big]^\top\z(t)\\
\ns\ds\qq\qq\q+Pb(t)+\big[C-D(R+D^\top PD)^\dag(B^\top P+D^\top PC+S)\big]^\top P\si(t)\\
\ns\ds\qq\qq\q+q(t)-\big[(R+D^\top PD)^\dag(B^\top P+D^\top PC+S)\big]^\top\rho(t)\Big\}dt+\z(t)dW(t),\q t\ges0,\ea\ee
such that
\bel{4.11}B^\top\eta(t)+D^\top\z(t)+D^\top P\si(t)+\rho(t)\in\sR(R+D^\top PD),\qq\ae t\in[0,\i),~\as\ee
In the above case, any closed-loop optimal strategy is given by the following:
\bel{optimal*}\left\{\2n\ba{ll}
\ds\bar\Th=-(R+D^\top PD)^\dag(B^\top P+D^\top PC+S)+\big[I-(R+D^\top PD)^\dag(R+D^\top PD)\big]\L,\\
\ns\ds\bar v(\cd)=-(R+D^\top PD)^\dag\big[B^\top\eta(\cd)+D^\top\z(\cd)+D^\top P\si(\cd)+\rho(\cd)\big]\\
\ns\ds\qq\qq+\big[I-(R+D^\top PD)^\dag(R+D^\top PD)\big]\n(\cd),\ea\right.\ee
for some $\L\in\dbR^{m\times n}$ and $\n(\cd)\in L^2_\dbF(0,\i;\dbR^n)$.

\et

\ms

Note that \rf{ARE-1} is the same as \rf{ARE1}, which means that the major condition under which Problem (EC) is solvable is the same as that for the closed-loop solvability of Problem (LQ)$_\i$. Moreover, with the expression $\bar\Th$ given in \rf{optimal*}, BSDE \rf{BSDE1} can be written as
\bel{BSDE2}\ba{ll}
\ns\ds d\eta(t)=-\Big\{(A+B\bar\Th)^\top\eta(t)+(C+D\bar\Th)^\top\z(t)+(C+D\bar\Th)^\top P\si(t)+Pb(t)+q(t)+\bar\Th^\top\rho(t)\Big\}dt\\
\ns\ds\qq\qq\qq+\z(t)dW(t),\qq t\ges0,\ea\ee
Hence, formally, if in the case that $b(\cd),\si(\cd),q(\cd)$ and $\rho(\cd)$ are all constants and the above BSDE has a constant solution $(\eta,0)$, then one has
$$(A+B\bar\Th)^\top\eta+(C+D\bar\Th)^\top P\si+Pb+q+\bar\Th^\top\rho=0,$$
which coincides with equation \rf{eta}. With the above, we also have
$$\bar v=-(R+D^\top PD)^\dag(B^\top\eta+D^\top P\si+\rho)+\big[I-(R+D^\top PD)^\dag(R+D^\top PD)\big]\n,$$
which coincides with the expression for $\bar v_0$ in \rf{v_0}. The above formal comparison shows that Problems (EC) and (LQ)$_\i$ are intrinsically match. Of course, this formal comparison does not lead to a rigorous ``equivalence'' proof between two problems.

\ms

Comparing Theorem \ref{finite-solvable} (ii) with Theorem \ref{Sun-Yong}, one can check the major solvability condition (the solvability of the algebraic Riccati equation) of Problem (EC) by checking the same condition for Problem (LQ)$_\infty$, which is easier since the latter is relatively easier to handle than the former. We point out that Theorem \ref{Sun-Yong} only gives an equivalence between the closed-loop solvability of Problem (LQ)$_\infty$ and the solvability of the algebraic Riccati equation \rf{ARE1} such that \rf{in Th} holds and BSDE \rf{BSDE1} has an adapted solution satisfying \rf{4.11}. It does not give condition(s) under which such a set of conditions can be fulfilled. In particular, it does not provide any useful information on conditions guaranteeing the solvability of the algebraic Riccati equation and strictly beyond the classical positive semi-definiteness condition \rf{Q>0}. We now, therefore, would like to take a closer look at sufficient conditions that ensure the solvability of the algebraic Riccati equation. We emphasize that the conditions should be strictly beyond the classical positive semi-definite condition. In another word, we want to include situations that either $Q$ or $R$ is negative definite.

\ms

We know that Problem (LQ)$_\i$ is a minimization problem of a quadratic functional on some Hilbert space. Thus the most natural sufficient condition is the uniform convexity of the functional. Clearly, such a condition is nothing to do with the nonhomogeneous terms $b(\cd),\si(\cd),q(\cd),\rho(\cd)$, and the initial state $x$. Hence, we need only to consider the homogeneous state equation with zero initial condition: (denoting $A_\Th=A+B\Th$ and $C_\Th=C+D\Th$)
\bel{state-3}\left\{\2n\ba{ll}
\ds dX_0^\Th(t)\1n=\1n\big[A_\Th X_0^\Th(t)\1n+\1n Bv(t)\big]dt\1n+\1n\big[C_\Th X_0^\Th(t)\1n+\1n Dv(t)\big]dW(t),\q t\ges0,\\
\ns\ds X_0^\Th(0)=0,\ea\right.\ee
and the cost functional
\bel{cost-3}\ba{ll}
\ns\ds\h J^{\,0,\Th}_\i(0;v(\cd))=\dbE\int_0^\i\[\lan(Q_\Th X_0^\Th(t),X_0^\Th(t)\ran
+2\lan S_\Th X_0^\Th(t),v(t)\ran+\lan Rv(t),v(t)\ran\]dt,\ea\ee
with
$$Q_\Th=Q+S^\top\Th+\Th^\top S+\Th^\top R\Th,\q S_\Th=S+R\Th.$$
The LQ problem associated the above is referred to as the {\it stabilized} LQ problem. Suppose $\Th\in\BTh[A,C;B,D]$ is so chosen that the following holds:
\bel{>d*}\h J^{\,0,\Th}_\i(0;v(\cd))\ges\d\dbE\int_0^\i|v(t)|^2dt,\qq\forall v(\cd)\in
\sU[0,\i),\ee
for some $\d>0$. Then from \cite{Sun-Yong 2018}, the corresponding homogeneous LQ problem admits a unique open-loop optimal control (for any initial condition), which is equivalent to the closed-loop solvability of the problem. Hence, the corresponding algebraic Riccati equation admits a solution $P\in\dbS^n$:
$$\ba{ll}
\ns\ds0=P(A+B\Th)+(A+B\Th)^\top P+(C+D\Th)^\top P(C+D\Th)+Q+S^\top\Th+\Th^\top S+\Th^\top R\Th\\
\ns\ds\qq-\big[B^\top P+D^\top P(C+D\Th)+S+R\Th\big]^\top(R+D^\top PD)^{-1}
\big[B^\top P+D^\top P(C+D\Th)+S+R\Th\big]\\
\ns\ds\q=PA+A^\top P+C^\top PC+Q+PB\Th+\Th^\top B^\top P+C^\top PD\Th+\Th^\top D^\top PC
+\Th^\top D^\top PD\Th\\
\ns\ds\qq+S^\top\Th+\Th^\top S+\Th^\top R\Th-(B^\top P+D^\top PC+S)^\top(R+D^\top PD)^{-1}(B^\top P+D^\top PC+S)\\
\ns\ds\qq-\Th^\top(B^\top P+D^\top PC+S)-(B^\top P+D^\top PC+S)^\top\Th-\Th^\top(R+D^\top PD)\Th\\
\ns\ds\q=PA+A^\top P+C^\top PC+Q-(B^\top P+D^\top PC+S)^\top(R+D^\top PD)^{-1}(B^\top P+D^\top PC+S),\ea$$
with
$$R+D^\top PD\ges\d I,$$
and
$$\ba{ll}
\ns\ds-(R+D^\top PD)^{-1}\big[B^\top P+D^\top P(C+D\Th)+S+R\Th\big]\\
\ns\ds=-(R+D^\top PD)^{-1}(B^\top P+D^\top PC+S)-\Th\in\BTh[A+B\Th,C+D\Th;B,D],\ea$$
which means that
$$-(R+D^\top PD)^{-1}(B^\top P+D^\top PC+S)\in\BTh[A,C;B,D].$$
Hence, in this case, conditions of (ii) in Theorem \ref{finite-solvable} hold and therefore Problem (EC) is solvable.

\ms

Now, we come to the question: When condition \rf{>d*} can be verified by an easily verifiable assumption, without assuming the classical positive-definiteness condition \rf{QR>0}? The following gives a sufficient condition which is inspired by a result from \cite{Sun-Xiong-Yong 2020}.

\bl{} \sl Let {\rm(H1)} hold and let $\Th\in\BTh[A,C;B,D]$, $Q_0\in\dbS^n$ with $Q_0>0$ such that the solution $\Pi$ to the following {\it Lyapunov inequality}:
\bel{Lyapunov}\Pi(A+B\Th)+(A+B\Th)^\top\Pi+(C+D\Th)^\top\Pi(C+D\Th)+S^\top\Th+\Th^\top S+\Th^\top R\Th+Q-Q_0\ges0,\ee
satisfies
\bel{R+DPD>0*}R+D^\top\Pi D -[\Pi B+(C+D\Th)^\top\Pi D+S^\top+\Th^\top R]^\top Q_0^{-1}
[\Pi B+(C+D\Th)^\top\Pi D+S^\top+\Th^\top R]\ges\d I,\ee
for some $\d>0$. Then \rf{>d*} holds.

\el

\begin{proof} \rm Let $X_0(\cd)\equiv X^\Th_0(\cd\,;0,v(\cd))$ be the solution to \rf{state-3}. Let $\Pi\in\dbS^n$ such that \rf{Lyapunov}--\rf{R+DPD>0*} hold. Observe the following:
$$\ba{ll}
\ns\ds d\lan\Pi X_0(t),X_0(t)\ran=\[\lan\Pi\((A+B\Th)X_0(t)+Bv(t)\),X_0(t)\ran
+\lan\Pi X_0(t),(A+B\Th)X_0(t)+Bv(t)\ran\\
\ns\ds\qq\qq\qq\qq\qq+\lan\Pi\([C+D\Th]X_0(t)+Dv(t)\),[C+D\Th]X_0(t)+Dv(t)\ran\]dt+[\cds]dW(t)\\
\ns\ds=\[\Big\langle\(\Pi(A+B\Th)+(A+B\Th)^\top\Pi+(C+D\Th)^\top\Pi(C+D\Th)\)X_0(t),X_0(t)\Big\rangle\\
\ns\ds\qq+2\big\langle\big[B^\top\Pi+D^\top\Pi(C+D\Th)\big]X_0(t),v(t)\big\rangle+\lan D^\top\Pi D
v(t),v(t)\ran\]dt+[\cds]dW(t).\ea$$
Then
$$\ba{ll}
\ns\ds\h J_\i^{\,0,\Th}(0;v(\cd))=\dbE\int_0^\infty\Big\{\Big\langle\[Q+S^\top\Th+\Th^\top S+\Th^\top R\Th+\Pi(A+B\Th)+(A+B\Th)^\top\Pi\\
\ns\ds\qq\qq\qq\qq\qq+(C+D\Th)^\top\Pi(C+D\Th)-Q_0\]X_0(t),X_0(t)\Big\rangle+\lan Q_0X_0(t),X_0(t)\ran\\
\ns\ds\qq\qq\qq\qq\qq+2\big\langle\1n\big[S\1n+\1n R\Th\1n+\1n B^\top\Pi\1n+\1n D^\top\Pi(C\1n+\1n D\Th)\big]X_0(t),v(t)\big\rangle
\1n+\1n\lan(R\1n+\1n D^\top\Pi D)v(t),v(t)\ran\1n\Big\}dt\\
\ns\ds\ges\dbE\int_0^\infty\Big\{|Q_0^{1\over2}X_0(t)|^2+2\big\langle\1n\big[S\1n+\1n R\Th\1n+\1n B^\top\Pi\1n+\1n D^\top\Pi(C\1n+\1n D\Th)\big]X_0(t),
v(t)\big\rangle\1n+\1n\lan(R\1n+\1n D^\top\Pi D)v(t),v(t)\ran\1n\Big\}dt\\
\ns\ds=\dbE\int_0^\infty \Big\{\big|Q_0^{1\over2}X_0(t)+Q_0^{-{1\over2}}\big[S^\top+\Th^\top R+\Pi^\top B
+(C^\top+\Th^\top D^\top)\Pi^\top D\big]v(t)\big|^2\\
\ns\ds\qq\qq\qq+\Big\langle\(R+D^\top\Pi D-[S+R\Th+B^\top\Pi+D^\top\Pi(C+D\Th)]Q_0^{-1}\\
\ns\ds\qq\qq\qq\qq\cd\big[S^\top+\Th^\top R+\Pi^\top B
+(C^\top+\Th^\top D^\top)\Pi^\top D\big]\)v(t),v(t)\Big\rangle\Big\}dt\ges\d\dbE\int_0^\infty|v(t)|^2dt.\ea$$
This proves \rf{>d*}.
\end{proof}

The above result gives some compatibility conditions among the coefficients of the state equation and the weighting matrices in the cost functional that ensure the uniform convexity condition \rf{>d*}. Let us take a closer look at \rf{Lyapunov} and \rf{R+DPD>0*}. Let us assume $R<0$ (or $R\ges0$ fails). Since $\BTh[A,C;B,D]\ne\varnothing$, we may find a $\Th\in\BTh[A,C;B,D]$ so that $[A+B\Th,C+D\Th]$ is stable. Hence, one can find a $\Pi>0$ so that
\bel{<0}\Pi(A+B\Th)+(A+B\Th)^\top\Pi+(C+D\Th)^\top\Pi(C+D\Th)<0.\ee
Note that the choices of $\Th,\Pi$ are independent of the weighting matrices $Q,S,R$ of the cost functional. Therefore, under the condition
\bel{R in D}\sR(R)\subseteq\sR(D),\ee
taking into account \rf{<0}, if necessary, replacing $\Pi$ by $\l\Pi$ for $\l>0$ large, we may have the following:
\bel{R+DPD>>0}R+D^\top\Pi D\ges2\d I,\ee
for some $\d>0$. Then we can find a large $Q_0>0$ such that \rf{R+DPD>0*} holds. Having
the $\Th,\Pi,Q_0$ given, we see that if $Q>0$ is sufficiently positive, the \rf{Lyapunov} will be true. This very rough analysis shows that $R<0$ could be compensated by the sufficient positiveness of $Q$ and the condition \rf{R in D}. Unfortunately, the above argument does not apply to the case $Q<0$. However, when $Q<0$, one should
expect certain compensation from the sufficient positiveness $R$. To see this, let us
recall \rf{ARE-0} which is rewritten here:
\bel{3.22*}\left\{\2n\ba{ll}
\ds\Pi_0A\1n+\1n A^\top\1n\Pi_0\1n+\1n C^\top\1n\Pi_0C\1n+\1n Q\1n-\1n(B^\top\1n\Pi_0
\1n+\1n D^\top\1n\Pi_0C\1n+\1n S)^\top\1n(R\1n+\1n D^\top\Pi_0D)^\dag\1n(B^\top\Pi_0\1n+\1n D^\top\1n\Pi_0C\1n+\1nS)\1n\ges\1n0,\\
\ns\ds R+D^\top\Pi_0D\ges0,\qq\sR(B^\top\Pi_0+D^\top\Pi_0C+S)\subseteq\sR(R+D^\top\Pi_0D).
\ea\right.\ee
Again, we look at the following two interesting cases:

\ms

\it Case 1. \rm Let $R<0$ (or $R\ges0$ fails). Pick a $\Pi_0\in\dbS^n$ with $\Pi_0>0$ so that
\bel{R+DPD>>0}R+D^\top\Pi_0 D\ges\d I,\ee
for some $\d>0$. For this, we still need \rf{R in D}. With such a $\Pi_0$, if $Q>0$ is sufficiently positive, then the following will hold:
\bel{}\Pi_0A\1n+\1n A^\top\1n\Pi_0\1n+\1n C^\top\1n\Pi_0C\1n+\1n Q\1n-\1n(B^\top\1n\Pi_0
\1n+\1n D^\top\1n\Pi_0C\1n+\1n S)^\top\1n(R\1n+\1n D^\top\Pi_0D)^{-1}(B^\top\Pi_0\1n+\1n D^\top\1n\Pi_0C\1n+\1nS)\1n\ges\d I,\ee
which is the form of \rf{3.22*} under condition \rf{R+DPD>>0}. This means that when $R<0$, as long as $Q$ is sufficiently positive, conditions of Theorem \ref{finite-solvable} (i)
are satisfied.

\ms

\it Case 2. \rm Let $Q<0$ (or $Q\ges0$ fails). If we can find a $\Pi_1\in\dbS^n$ such that
$$\Pi_1 A+A^\top\Pi_1+C^\top \Pi_1 C>0\qq\text{ or }\qq\Pi_1 A+A^\top\Pi_1+C^\top \Pi_1 C<0, $$
then we can find an $\a$ ($\a$ is negative if in the above, the second inequality holds)
and $R>0$ such that the first inequality in \rf{3.22*} holds for $\Pi_0=\a\Pi_1$. Next,
if $R$ is positive enough, we will have the second inequality in \rf{3.22*}.

\ms

From those observations, we see that conditions of Theorem \ref{finite-solvable} (i)
can be verified even if one of $Q$ and $R$ is negative definite. Also, the above Case 2 suggests us that when $R>0$ is sufficiently positive, then we may take some $\d>0$ and check the condition \rf{3.22*} with $R$ replaced by $R-\d I$. If such a condition is
satisfied, then Problem (EC) will be solvable. We will present results relevant to this in the following section.

\section{Optimal Value and Regularization of Ergodic Problem}

In Theorem \ref{finite-solvable}, we have proved that under (H2), Problem (EC) is finite. In this section, we will find the optimal value $\sE$ in this case. First, let us refine Theorem \ref{finite-solvable} (i). To this end, we introduce the following hypothesis which is a part of (H2).

\ms

{\bf(H2)$'$} Let (H1) and \rf{ARE-0} hold.

\ms

Let us again look at the homogeneous problem associated with \rf{state-3}--\rf{cost-3}, for some $\Th\in\BTh[A,C;B,D]$. We have the following result.

\bp{Prop-5.1} \sl Let {\rm(H2)$'$} hold. Then for any $(\Th,v(\cd))\in\dbU[0,\i)$,
\bel{J>0}\h J^{\,0,\Th}_\infty(0;v(\cd))\ges0.\ee

\ep

\begin{proof} Let $P\in\dbS^n$ and we apply It\^o's formula to $\lan PX^\Th_0(\cd),X^\Th_0(\cd)\ran$.
$$\ba{ll}
\ns\ds0=\dbE\int_0^\i\[\lan P[A_\Th X^\Th_0(t)+Bv(t)],X^\Th_0(t)\ran
+\lan PX^\Th_0(t),A_\Th X^\Th_0(t)+Bv(t)\ran\\
\ns\ds\qq\qq\qq +\lan P[C_\Th X^\Th_0(t)+Dv(t)],C_\Th X^\Th_0(t)+Dv(t)\ran\]dt\\
\ns\ds\q=\dbE\int_0^\i\[\Big\langle\(P A_\Th+A_\Th^\top P+C_\Th^\top PC_\Th\)X^\Th_0(t),X^\Th_0(t)\Big\rangle\\
\ns\ds\qq\qq\qq +2\Big\langle\(B^\top P+D^\top PC_\Th\)X^\Th_0(t),v(t)\Big\rangle+\lan D^\top PDv(t),v(t)\ran\]dt.\ea$$
Hence,
$$\ba{ll}
\ns\ds\h J^{\,0,\Th}_\i(0;v(\cd))=\dbE\int_0^\i\[\lan(Q_\Th+PA_\Th+A_\Th^\top P+C_\Th^\top PC_\Th)X^\Th_0(t),X^\Th_0(t)\ran\\
\ns\ds\qq\qq\qq\qq\qq+2\lan(S_\Th+B^\top P+D^\top PC_\Th)X^\Th_0(t),v(t)\ran+\lan(R+D^\top PD)v(t),v(t)\ran\]dt.\ea$$
Let $P\in\dbS^n$ satisfy \rf{ARE-0}. Then
$$\sR\big(S_\Th+B^\top P+D^\top PC_\Th\big)=\sR\big(S+B^\top P+D^\top PC+(R+D^\top PD)\Th\big)\subseteq\sR(R+D^\top PD).$$
Hence, we may complete the square to obtain
$$\ba{ll}
\ns\ds\h J^{\,0,\Th}_\i(0;v(\cd))=\dbE\int_0^\i\[ \Big\langle\(Q_\Th+PA_\Th+A_\Th^\top P+C_\Th^\top PC_\Th-(S_\Th+B^\top P+D^\top PC_\Th)^\top(R+D^\top PD)^\dag\\
\ns\ds\qq\qq\qq\qq\qq\q\cd(S_\Th+B^\top P+D^\top PC_\Th)\)X^\Th_0(t),X^\Th_0(t)\Big\rangle\\
\ns\ds\qq\qq\qq\qq\qq+\Big|(R+D^\top PD)^{1\over2}\(v(t)+(R+D^\top PD)^\dag(S_\Th+B^\top P+D^\top PC_\Th)X^\Th_0(t)\)\Big|^2\]dt.\ea$$
Note that
$$\ba{ll}
\ns\ds Q_\Th+PA_\Th+A_\Th^\top P+C_\Th^\top PC_\Th-(S_\Th+B^\top P+D^\top PC_\Th)^\top(R+D^\top PD)^\dag(S_\Th+B^\top P+D^\top PC_\Th)\\
\ns\ds=Q+S^\top\Th+\Th^\top S+\Th^\top R\Th+P(A+B\Th)+(A+B\Th)^\top P+(C+D\Th)^\top P(C+D\Th)\\
\ns\ds\q-[S+R\Th+B^\top P+D^\top P(C+D\Th)]^\top(R+D^\top PD)^\dag[S+R\Th+B^\top P+D^\top P(C+D\Th)]\\
\ns\ds=Q+PA+A^\top P+C^\top PC+(S^\top+PB+C^\top PD)\Th+\Th^\top(S+B^\top+D^\top PC)
+\Th^\top(R+D^\top PD)\Th\\
\ns\ds\q-(S+B^\top P+D^\top PC)^\top(R+D^\top PD)^\dag(S+B^\top P+D^\top PC)\\
\ns\ds\q-\Th^\top(R+D^\top PD)(R+D^\top PD)^\dag(S+B^\top P+D^\top C)\\
\ns\ds\q-(S+B^\top P+D^\top PC)^\top(R+D^\top PD)^\dag(R+D^\top PD)\Th-\Th^\top(R+D^\top PD)\Th\\
\ns\ds=Q+PA+A^\top P+C^\top PC-(S+B^\top P+D^\top PC)^\top(R+D^\top PD)^\dag(S+B^\top P+D^\top PC)\ges0.\ea$$
Hence, \rf{J>0} follows. \end{proof}

\ms

Now, for any $\d>0$, we denote
$$R_\d=R+\d I,$$
and
$$g_\d(x,u)=g(x,u)+\d|u|^2\equiv\lan Qx,x\ran+2\lan Sx,u\ran+\lan R_\d u,u\ran+2\lan q,x\ran+2\lan\rho,u\ran.$$
Then, correspondingly, we introduce the following {\it regularized} ergodic cost functional
\bel{cE-d}\cE_\d(\Th,v)=\int_{\dbR^n}\(g(x,\Th x+v)+\d|\Th x+v|^2\)\pi^{\Th,v}(dx),\ee
and introduce the following optimal control problem.

\ms

{\bf Problem (EC)$_\d$.} Let (H1) hold. Find a $(\bar\Th,\bar v)\in\dbU$ such that
$$\cE_\d(\bar\Th,\bar v)=\inf_{(\Th,v)\in\dbU}\cE_\d(\Th,v)\equiv\sE_\d.$$

\ms

The following lemma reveals the relationship between Problems (EC) and (EC)$_\d$.

\bl{valueconvergence} \sl Let {\rm(H1)} hold. If $\sE$ is finite, then
\bel{seconvergence}\sE=\lim_{\d\rightarrow 0^+}\sE_\d,\ee
which is true if {\rm(H2)} holds. In particular, this is true if $g$ is bounded below.
\el

\begin{proof} It is easy to see that $\sE\les\sE_\d$. Next, let $(\Th_k,v_k)\in\dbU$ be a minimizing sequence of $\cE(\cd)$ such that
$$\sE\les\cE(\Th_k,v_k)\equiv\int_{\dbR^n}g(x,\Th_kx+v_k)\pi_k(dx)<\sE+{1\over k},\qq k\ges1,$$
where $\pi_k$ is the invariant measure corresponding to $(\Th_k,v_k)\in\dbU$. Then
\bel{lowerboundse-2}\ba{ll}
\ns\ds\sE\ad\ges\cE(\Th_k,v_k)-{1\over k}\equiv\int_{\dbR^n}g(x,\Th_kx+v_k)\pi_k(dx)-{1\over k}\\
\ns\ad=\limsup_{\d\to0+}\(\sE_\d-\d \int_{\dbR^n}\big|\Th_kx+v_k|^2\pi_k(dx)\)-{1\over k}
=\limsup_{\d\to0^+}\sE_\d-{1\over k}.\ea\ee
Here we have used the fact that $\int_{\dbR^n}|x|^2\pi_k(dx)<\infty$. Then, one has
$$\sE\ges\limsup_{\d\to0+}\sE_\d\ges\sE,$$
proving \rf{seconvergence}.  \end{proof}

Let us call the LQ problem with $R$ replaced by $R_\d$ Problem (LQ)$_{\i,\d}$. The cost functional of this problem reads
$$\h J_{\i,\d}(x;u(\cd))=\h J_\i(x;u(\cd))+\d\dbE\int_0^\i|u(t)|^2dt,\qq\forall u(\cd)\in\sU_{ad}[0,\i).$$
Hence, by Proposition \ref{Prop-5.1}, the cost functional of the corresponding homogeneous problem satisfies the following:
$$\h J^{\,0,\Th}_{\i,\d}(0;v(\cd))=\h J^{\,0,\Th}_\i(0;v(\cd))+\d\dbE\int_0^\i|v(t)|^2dt\ges\d\dbE\int_0^\i|v(t)|^2dt,\q\forall
(\Th,v(\cd))\in\dbU[0,\i).$$
Hence, Problem (LQ)$_{\i,\d}$ is uniquely closed-loop solvable. Consequently, the following algebraic Riccati equation admits a solution $\h P_\d$:
\bel{ARE-d}\ba{ll}
\ns\ds\h P_\d A+A^\top\h P_\d+C^\top\h P_\d C+Q \\
\ns\ds\qq\qq\q-(B^\top\h P_\d+D^\top\h P_\d C+S)^\top(R_\d+D^\top\h P_\d D)^{-1}(B^\top \h P_\d+D^\top\h P_\d C+S)=0,\ea\ee
with
\bel{Stabilizer}\h\Th_\d=-(R_\d+D^\top\h P_\d D)^{-1}(B^\top\h P_\d+D^\top\h P_\d C+S)\in\BTh[A,C;B,D].\ee
Then applying Theorem \ref{finite-solvable} (ii), we have the optimal strategy of Problem (EC)$_\d$ given by the following:
\bel{optd}\left\{\2n\ba{ll}
\ds\h\Th_\d=-(R_\d+D^\top\h P_\d D)^{-1}(B^\top\h P_\d+D^\top\h P_\d C+S),\\
\ns\ds\h v_\d=-(R_\d+D^\top\h P_\d D)^{-1}(B^\top\h\eta_\d+D^\top\h P_\d\si+\rho),\\
\ns\ds\h\eta_\d=-\big[(A+B\h\Th_\d)^\top\big]^{-1}\big[\h P_\d b+(C+D\h\Th_\d)^\top
\h P_\d\si+q+\h\Th_\d\rho\big],\ea\right.\ee
with the optimal value:
$$\sE_\d\equiv\cE_\d (\h\Th_\d,\h v_\d)=\lan\h P_\d\si,\si\ran+2\lan\h\eta_\d,b\ran-\lan(R_\d+D^\top\h P_\d D)\h v_\d,\h v_\d\ran.$$

\ms

Now we can present the approximation theorem for the value of  Problem (EC).

\begin{theorem}\label{mainth}\sl Let {\rm(H2)} hold. Then Problem (EC) is finite with
\bel{covergencesed}\sE=\lim_{\d\to0^+}\cE(\h\Th_\d,\h v_\d)=\lim_{\d\to0^+}\cE_\d(\h\Th_\d,\h v_\d)=\lim_{\d\to0^+}\[\lan\h P_\d\si,\si\ran+2\lan\h\eta_\d,b\ran-\lan(R_\d+D^\top\h P_\d D) \h v_\d,\h v_\d\ran\].\ee
i.e., $(\h\Th_\d,\h v_\d)\in\dbU$ is a minimizing sequence of Problem (EC). Moreover, if  $(\h\Th_\d,\h v_\d)$ has a convergent subsequence with limit $(\h\Th,\h v)\in\dbU$, then Problem (EC) is solvable and $(\h\Th,\h v)\in\dbU$ is an optimal strategy.
	
\end{theorem}

The proof is clear.

\ms

To conclude this section, we point out the steps to obtain the optimal value of Problem (EC) as follows:

\ms

{\it Step-1.} Consider homogeneous infinite-horizon optimization problem
$$dX(t)=\big[AX(t)+Bu(t)\big]dt+\big[CX(t)+Du(t)\big]dW(t),$$
with cost functional
$$\h J_{\i,\d}(x;u(\cd))=\dbE\int_0^\i\big[\lan QX(t),X(t)\ran+2\lan SX(t),u(t)\ran+\lan R_\d u(t),u(t)\ran\big]dt.$$
By \cite{Sun-Yong 2018}, we know that the closed-loop and open-loop solvability are equivalent. Given (H2), such problem is solvable and we can find an optimal strategy $(\h\Th_\d,0)\in\dbU$.

\ms

{\it Step-2.} Consider the non-homogeneous optimization problem with respect to $v$ in the problem with the state equation:
$$dX(t)=\big[(A+B\h\Th_\d)X(t)+Bv+b\big]dt+\big[(C+D\h\Th_\d)X(t)+Dv+\si\big]dW(t),$$
and the cost function $v\mapsto\cE_\d(\h\Th_\d,v)$, i.e., find an optimal $\h v_\d\in\dbR^m$ such that
$$\cE_\d(\h\Th_\d,\h v_\d )=\inf_{v\in\dbR^m}\cE_\d(\h\Th_\d,v).$$

{\it Step-3.} Let $\d\rightarrow 0^+$ and obtain $\sE=\lim\limits_{\d\to0^+}\cE_\d(\h\Th_\d,\h v_\d )$.

\ms

Note that we cannot apply such method to Problem (EC) directly, since the first step is not necessarily going through under (H2) if $R_\d$ is replaced by $R$.

\section{Examples}\label{sec:exp}

In the section, we will present two one-dimensional examples to illustrate our results.\ss

{\bf Example 1}. \rm  Consider the following one-dimensional controlled SDE:
$$dX(t)=\big[AX(t)+u(t)+b\big]dt+\big[CX(t)+\si\big]dW(t),$$
with cost functional rate
$$g(x,u)=Q x^2+2Sxu.$$
Note that we assume $B=1$, $D=R=\rho=q=0$ in such example.  We can see that $R+D^\top PD=0$ and the classical algebraic Riccati equation does not hold. Through a direct calculation which will be presented in the Appendix, we list all the possible cases in the Table 1 below, where
\bel{v_Th} v_{\Th}\equiv-\frac{(2(A+\Th)+C^2)}{2[Q-S(2A+C^2)]}\(Sb+\frac{(Q+2S\Th)C\si}{2(A+\Th)+C^2}\).\ee

	\begin{table}[!htb]

	\centering
	
	\begin{tabular}{|c|c|c|c|c|}
		\hline
		~ & { $S(2A+C^2)$}&$CSb+(Q-2AS)\si$ &{\bf Finite}&{\bf Solvable} \cr\cline{3-5}
		\hline
		\multirow{2}{*}{I} &\multirow{2}{*}{ $<Q$}& $=0$&Yes&$(\Theta,v_{\Th})$ is optimal for any $\Th\in\BTh[A,C;B,D]$ \cr
		\cline{3-5}
		&~&$\neq0$ &Yes&No\cr\hline
		\multirow{2}{*}{II} &\multirow{2}{*}{ $=Q$}& $=0$&Yes&any $(\Th,v)\in\BTh[A,C;B,D]\times\dbR^m$ is optimal\cr
		\cline{3-5}
		&~&$\neq0$ &No&No\cr\hline
		{III} & { $>Q$}& &No&No\\
		\cline{3-5}
		\hline
	\end{tabular}
\caption{Finiteness and Solvability of Problem (EC) with $R=D=0$.}
\end{table}
%

%
	%
	%Note that $$h(-\infty)=-\frac{S^2(b+C\sigma)^2}{Q-2AS-SC^2}-S\sigma^2.$$
	%If there exists a $-\infty<\Theta_0<-A-C^2/2$ such  that $h(\Theta_0)\leq h(-\infty)$ the problem is  solvable.
	%Otherwise, the problem is not  solvable. For example, if $C\neq0$ and %$\sigma=0$,  the problem is not solvable. If $C=0$ and $\sigma\neq 0$,
	%$h(\Theta)$ is an increasing function and thus the problem is not  solvable.
	% If $C=0$ and $\sigma=0$, the problem is solvable with any $\Theta<-A$.\ms	

We will see that our assumption (H2) corresponds to  case I and $CSb+(Q-2AS)\si=CS(b+C\si)=0$ in case II. Observed
from case I, we can see that $S(2A+C^2)<Q$ is sufficient for finiteness. Therefore, if $S(2A+C^2)<0$ is negative, $Q$ is allowed to be negative, even if $R=0$. This is the case, if the system $[A,C]$ is stable (which implies $2A+C^2<0$) and $S>0$. In such a case, $S(2A+C^2)$ gives a lower bound for $Q$ so that Problem (EC) is finite. It is more  surprising that as long as $2A+C^2\ne0$ with $S$ having the opposite sign, then $S(2A+C^2)<0$ and which allows $Q$ to be negative. Hence, $S$ has a contribution to the finiteness of Problem (EC). For example, if $A=B=C=1$, $D=0$, then we may allow
$$\begin{pmatrix}Q&S^\top\\ S&R\end{pmatrix}=\begin{pmatrix}-1&-1\\-1&0\end{pmatrix},\qq
g(x,u)=-x^2-2xu.$$
Therefore, the function $g(x,u)$ is even unbounded below. This shows that our assumptions assumed in Theorem \ref{finite-solvable} are much weaker than one could imagine.

\ss

 %Moreover , we can see that if $S\ne0$,
%
%$$\begin{pmatrix}Q&S^\top\\S&R\end{pmatrix}=\begin{pmatrix}Q&S^\top\\S&0\end{pmatrix}$$ %
%is indefinite.
Now we try to use  Theorem \ref{finite-solvable} conclude the finiteness and solvability of Problem (EC). We can compare our results with Table 1.\ss

{\noindent\it  Case I:}  $S(2A+C^2)<Q$. Observed from   (H2), we take a $\Pi_0=-S$, then
$$\ba{ll}\left\{\2n\ba{ll}
\ds (2A+C^2)\Pi_0+Q>0;\\
\ns\ds\sR(\Pi_0+S )\subseteq\sR(0); \\
\ns\ds 0\ges0.\ea\right.\ea$$
One can see that $(\Pi_0,\eta_0)=(-S,0)$ can guarantee (H2). Thus our theorem says that  if $R=0$ and $Q>S(2A+C^2)$, $\sE$ is finite. \ss

\ms

{\noindent\it  Case II:} $S(2A+C^2)=Q$. Note that (H3) requires us to take $\Pi_0=-S$, $\eta_0=0$ and
$C\sigma+b=0.$ Thus  our theorem states that Problem (EC) is solvable if $Q=S(2A+C^2)$ and
$C\sigma+b=0.$ This corresponds to $CSb+(Q-2AS)\si=0$ in case II from Table 1.  We also can see that such case is a special class of (H2).\ss

 While we can see that Problem (EC) is solvable in the case  $CSb+(Q-2AS)\si=0$ in case I  where (H3) is not fulfilled. Therefore (H3) is not necessary for solvability of Problem (EC).\ss

\ms

Now we will present how the regularized problem approximates Problem (EC). Firstly let us assume $S(2A+C^2)<Q$. As stated previously,  (H2) is verified by $\Pi_0=-S$ and $\eta_0=0$.
		Write $$\a=Q-S(2A+C^2)\text{ and }\b ={2A+C^2\over2}.$$
	The algebraic Riccati equation  writes
	$$ (2A+C^2)P+Q-\d^{-1}(P+S)^2=0,$$
	which admits two solutions
	$$P=-S+\d \b\pm\sqrt{\d\a+\d^2\b ^2}.$$	
Then
	$$\Theta=-\d^{-1}(-S+\d\b\pm\sqrt{\d\a+\d^2\b ^2}+S)=-\b \pm \sqrt{\d^{-1}\a+\b^2}.$$
	We have to select  $\Theta$ to  stabilize the system, i.e.
	$$\h\Theta_\d=-\b - \sqrt{\d^{-1}\a+\b ^2}\text{ and }
	\h P_\d=-S+\d\b+\sqrt{\d\a+\d^2\b^2}.$$
	Then it follows that
	$$\left\{\ba{ll}\ns\ad\h\eta_\d=-(A-\b -\sqrt{\d^{-1}\a+\b ^2})^{-1}(-S+\d\b+\sqrt{\d\a+\d^2\b ^2})(b+C\si);\\ [1mm]
	\ns\ad \h v_\d=\d^{-1}(A-\b -\sqrt{\d^{-1}\a+\b ^2})^{-1}(-S+\d\b+\sqrt{\d\a+\d^2\b^2}) (b+C\si).\ea\right.$$
	As a result, as $\d\to0^+,$
	$$\ba{ll}\sE_\d\ad= \h P_\d\si^2+2b\h\eta_\d-\d\h v_\d^2\\
	\ns\ad=-(S+\d \b+\sqrt{\d\a+\d^2\b ^2})\si^2-2b(A-\b -b\sqrt{\d^{-1}\a+\b^2 })^{-1}(-S+\d\b+\sqrt{\d\a+\d^2\b ^2})(b+C\si)\\
	\ns\ad\q -\d^{-1}\[(A-\b - \sqrt{\d^{-1}\a+\b^2 })^{-1}(-S+\d\b+\sqrt{\d\a+\d^2\b ^2})(b+C\si)\]^2\\
	\ns\ad\rightarrow -\frac{S^2(b+C\si)^2}{Q-S(2A+C^2)}-S\si^2=\sE.\ea$$
	This verifies the approximation procedure in Theorem \ref{mainth}.\ss
	
	 From Table 1, we know that when $CSb+(Q-2AS)\si=0$ (e.g. $b=\si=0$), Problem (EC) is solvable. %One can see that the case $A=C=S=\si=1$, $b=-2$ and $Q=4$ satisfies 2-V.
	 We can see that since $\a>0$,  $\h\Theta_\d$  blows up with a  rate of $\sqrt{\d^{-1}}.$   Hence the convergence of $(\h\Th_\d,\h v_\d)$ is unnecessary for solvability of Problem (EC) where the sufficiency is stated in Theorem \ref{mainth}.\ss
	
	In the case $\a=Q-S(2A+C^2)=0$ and $b+C\si=0$, one can see that $\Th$ and $v$ are bounded and the limit is an optimal strategy. This coincides with our results in Theorem \ref{mainth}.\ss

	If $Q-S(2A+C^2)=0$ and $b+C\si\neq 0$, $\sE_\d\to-\infty.$ The problem is not finite essentially.
	\ms
	
	 {\bf Example 2.} In this example, we will deal with the case $D\neq 0$. We consider one-dimensional SDE
	$$dX(t)=(AX(t)+Bu(t)+b)dt+(CX(t)+Du(t)+\sigma)dW(t)$$
	with  $$g(x,u)=Q x^2+2Sxu+Ru^2+2qx+2\rho u.$$
We use the following notations in this example only
$$\left\{\ba{ll}\a=D^{-2}(B+CD)^2-(2A+C^2);\\
\b=Q-D^{-2}(2A+C^2)R-2D^{-2}[S-D^{-2}R(B+CD)](B+CD);\\
\g=D^{-2}[D^{-2}R(B+CD)-S]^2.\ea\right.$$
We can see that (H1) is equivalent to $\a>0$ and $\Th\in\BTh[A,C;B,D]$ if and only if
$$|D^2\Th+(B+CD)|<\sqrt\a\,|D|.$$
The results for finiteness using (H2) and solvability using (H3) are presented in  Table 2
where
$$\Theta_*=-D^{-2}(B+CD)-|D|^{-1}\sqrt{\a}\cd\text{sgn}\big\{[S-D^{-2}R(B+CD)]\big\}.$$
Here we note that $\Th_*\notin\BTh[A,C;B,D]$.

\begin{table}[!htb]\label{table2}
	\begin{center}
		\begin{tabular}{|c|c|c|c|}

			\hline
			I& $\beta-2\sqrt{\a\g}>0$	 & Solvable\\\hline
			\multirow{2}{*}{II}&$\b=\g=0$ and
			 	& \multirow{2}{*}{Solvable}\\
			&   	$B\eta_0-D^{-1}R\si+\rho=A\eta_0+q-D^{-2}R(b+C\si)=0$  for some $\eta_0$          &\\
			\hline
			\multirow{2}{*}{III}&$\g\neq0,\q\b-2\sqrt{\a\g}=0$, and  &\multirow{2}{*}{Finite}\\
			 &$q+\Th_*\rho+[b+(C+D\Th_*)\si](\sqrt{\gamma/\a}-D^{-2}R)\in\sR(A+B\Th_*)$&\\
			\hline
		\end{tabular}\ss
		
	\end{center}
	\caption{Finiteness  and solvability of Problem (EC) with $D\neq 0$.}
\end{table}	

\no Note that
$$\ba{ll}
\ds\beta-2\sqrt{\a\g}\\
\ns\ds=\left\{\2n\ba{ll}\ds Q+D^{-4}R(D\sqrt\a-B-CD)^2+2D^{-2}S(D\sqrt\a-B-CD),~\text{if } D^{-1}(S-D^{-2}R(B+CD))\les0;\\
%\ns\ad Q-R(2A+C^2)=Q+R\a-S(B+C),\q\text{if } S-R(B+C)=0;\\
\ns\ds Q+D^{-4}R(D\sqrt\a+B+CD)^2-2D^{-2}S(D\sqrt\a+B+CD),~\text{if } D^{-1}(S-D^{-2}R(B+CD))>0. \ea\right.\ea$$
%Thus 2-I tells us that the problem
%is solvable when $Q,R$ are large.
%For the cases {\rm (ii), (iii),} and {\rm (iv)}, the unique optimal strategy is
%\bel{sta}\Th=\frac{2\a((B+C)R-S)}{\beta+\sqrt{\beta^2-2\a\gamma}}.\eel
We can see if the third term is positive, Problem (EC) can be solvable even if $Q$ and $R$ are all negative. For example, if $A=B=C=D=1$, then $\a=1$ and
$$\ba{ll}
\ns\ds D^{-1}\big[S-D^{-2}R(B+CD)\big]=S-2R<0,\\
\ns\ds\b-2\sqrt{\a\g}=Q+R(-1)^2+2S(-1)=Q+R-2S>0,\ea$$
provided, say, $Q=R=-1$ and $S=-{5\over2}$. According to the above, we have the solvability of the corresponding Problem (EC). Interestingly, in the current case, we have
$$\begin{pmatrix}Q&S^\top\\ S&R\end{pmatrix}=\begin{pmatrix}-1&-{5\over2}\\-{5\over2}&-1\end{pmatrix},\qq
g(x,u)=-x^2-5xu-u^2.$$
Both $Q$ and $R$ are negative!

\ms

\section{Concluding Remarks}\label{sec:con}

In the paper, we have explored the ergodic optimal control problems for linear systems with  quadratic costs. Compared to the previous works on similar problems, we deal with a class of ergodic control problems allowing the weighting matrices of the cost functional to be indefinite. We have presented sufficient conditions for finiteness and solvability of the ergodic control problem. Comparing Problem (EC) with classical LQ problem on $[0,\i)$, we see that the algebraic Riccati equation in the condition for the solvability of the former coincides with that for the latter. Further, we have found a general sufficient condition under which the quadratic cost functional of the stabilized LQ problem is uniformly convex, which will lead to the closed-loop solvability of the LQ problem and therefore the solvability of Problem (EC). Moreover, when the problem is merely finite, we find a way of finding the optimal value of the problem by solving a sequence of regularized ergodic problems. Examples for one-dimensional cases showed that conditions that we have found for the solvability of Problem (EC) are sufficient but not necessary (see the discussion on cases I and II in Example 1). The difficulty of finding equivalent conditions for the solvability of Problem (EC) is probably due to the fact that the ergodic cost function $\cE(\Th,v)$, involving invariant measure, is not a convex function of $(\Th,v)$. We hope to report on the investigation of this in the near future.

\section*{Appendix}
In the appendix, we will present the proof in our examples.
\ms

{\noindent \it Proof of Example 1.} Now let us present the proof for Table 1.	It is easy to see that $$\dbU=\Big\{u(x)=\Th x+v:\Th<-{2A+C^2\over2}\Big\}.$$
Take $u(x)=\Theta x+v\in \dbU$. Through the use of It\^o's formula, simple calculation yields the first and the second moment of the invariant measure are
$$m_1=-\frac{b+v}{A+\Theta},\q m_2=\frac{2(b+v)^2}{[2(A+\Th)+C^2](A+\Th)}-\frac{\sigma^2}{2(A+\Th)+C^2}
+\frac{2C\sigma(b+v)}{(A+\Th)[2(A+\Th)+C^2]}.$$
Then,
$$\cE(\Th,v)=2\frac{Q-S(2A+C^2)}{2(A+\Th)+C^2}\frac{(b+v)^2}{A+\Theta}+\frac{2(b+v)}
{(A+\Th)}\(Sb+\frac{(Q+2S\Theta)C\sigma}{2(A+\Th)+C^2}\)-\frac{(Q+2S\Theta)\sigma^2}
{2(A+\Th)+C^2}.$$
Recall that $A+\Theta<0$ and $2(A+\Th)+C^2<0$.\ss

(1) If $S(2A+C^2)>Q$, the coefficient of the quadratic term is negative.  Problem (EC) is not finite.\ss

(2) If $S(2A+C^2)=Q$, then $$\cE(\Theta,v)=\frac{2(b+v)}{(A+\Theta)}S(b+C\sigma)-\frac{(Q+2S\Th)\sigma^2}{2(A+\Th)+C^2}.$$
If  $b+C\sigma=0$,  the problem is  solvable  and any admissible strategy is optimal. If $b+C\sigma\neq 0$, the problem is not finite.\ms

(3) Assume $S(2A+C^2)<Q$. Since $v$ can be taken arbitrarily, it follows that
$$\ba{ll}
\ns\ds h(\Th)=\inf_v\cE(\Theta,v)= -\frac{2(A+\Th)+C^2}{2(A+\Theta)[Q -S(2A+C^2)]}\(Sb+\frac{(Q+2S\Theta)C\sigma}{2(A+\Th)+C^2}\)^2-\frac{(Q+2S\Theta)\sigma^2}
{2(A+\Th)+C^2}\\
\ns\ds=-\frac{(2A+2\Theta+C^2)}{2(A+\Theta)[Q -S(2A+C^2)]}\(S(b+C\si)+\frac{[Q-S(2AS+C^2)]C\sigma}{2(A+\Th)+C^2}\)^2
-\frac{[Q-S(2A-C^2)]\si^2}{2(A+\Th)+C^2}-S\si^2\\
\ns\ds=-\frac{[2(A+\Th)+C^2]S^2(b+C\si)^2}{2(A+\Theta)[Q-S(2A+C^2)]}-\frac{C\si S(b+C\si)}{A+\Theta}-\frac{[Q-S(2A+C^2)]\sigma^2}{2(A+\Th)+C^2}{2(A+\Th)+C^2\over2(A+\Th)}-S\si^2,\ea$$
where the minimum is taken at
$$v_\Th=-\frac{2(A+\Th)+C^2}{2[Q-S(2A+C^2)]}\(Sb+\frac{(Q+2S\Theta)C\sigma}{2(A+\Th)+C^2}\).$$
One can see $h(\Theta)$ is bounded from below.
Thus the problem is finite if $2AS+SC^2<Q$.\ss

Now we want to look at when it is solvable. Note that
$$h(-\infty)=-\frac{S^2(b+C\si)^2}{Q-S(2A+C^2)}-S\si^2.$$
We solve the inequality
$$h(\Theta)\les h(-\infty).$$
This is equivalent to
$$\ba{ll}\ns\ad-\frac{(2(A+\Th)+C^2)^2S^2(b+C\si)^2}{Q-S(2A+C^2)}-2C\si S(b+C\si)[2(A+\Th)+C^2]-[Q-S(2A+C^2)]\si^2[2(A+\Th)+C^2]\\
\ns\ad\les -2\frac{S^2(b+C\si)^2(A+\Th)[2(A+\Th)+C^2]}{Q-S(2A+C^2)}\\
\ns\ad=-\frac{S^2(b+C\si)^2[2(A+\Th)+C^2]^2}{Q-S(2A+C^2)}+\frac{C^2S^2(b+C\si)^2[2(A+\Th)
+C^2]}{Q-S(2A+C^2)}.\ea$$
Straightforward calculation yields that
$$[2(A+\Th)+C^2]\(\frac{C^2S^2(b+C\si)^2}{Q-S(2A+C^2)}+2C\si S(b+C\si)+[Q-S(2A+C^2)]\si^2\)\ges 0.$$	
Note that the second term is a perfect square and $2(A+\Th)+C^2<0$. It  holds for some $\Theta<A+ {C^2} /2$ if and only if the second term is 0, i.e.
$$CS(b+C\si)+[Q-S(2A+C^2)]\si=0.$$
In this case, any admissible strategy is optimal. Otherwise, the problem is finite but not solvable.\ms

{\noindent \it Proof of Example 2.}
(1) {\bf  Finiteness}.  (H2) is equivalent to that there exist  $\Theta\in\Upsilon[P]$ and $\eta_0 $ such that
	\bel{griccati-ass0000}\left\{\2n\ba{ll}
	\ds \Psi[P]=(2A+C^2)P+Q-(P(B+CD)+S )^2(R+D^2P)^{\dagger}\geq  0;\\
	\ns\ds\sR( P(B+CD)+S)\subseteq\sR(R+D^2P); \\
	\ns\ds R+ D^2P\ges0;\\
	%\ns\ds P\si+B\eta\in\sR(R+P)\\
	\ns\ds DP \sigma+B\eta_0+\rho\in\sR(R+D^2P);\\
	\ns\ds (A+B\Theta)\eta_0+q+P b+\Theta\rho+(C+D\Theta)P\si\in\sR(\Psi[P]).\ea\right.\ee
	Note that for $P>-D^{-2}R$,
	$$\ba{ll}\Psi[P]\ad=[(2A+C^2)-D^{-2}(B+CD)^2](P+D^{-2}R)-\frac{D^{-2}(S-D^{-2}R(B+CD))^2}{P+D^{-2}R}\\
	\ns\ad\q+Q-D^{-2}(2A+C^2)R-2D^{-2}(S-D^{-2}R(B+CD))(B+CD)\\
	\ns\ad=-\a(P+D^{-2}R)-\gamma(P+D^{-2}R)^{-1}+\beta.
	%\ns\ad \leq -2\sqrt{-(2A+C^2-(B+C)^2)(S+R(B+C))^2}
	\ea$$

	If $\g=0$,  $\Psi[P]$ is decreasing of $P$. Thus we need  $\b\geq0$.
	If $\b>0$, we can select a $P=-D^{-2}R+\d$ for small $\d$. Then \eqref{griccati-ass0000} holds.
	If $\b=0$, equation \eqref{griccati-ass0000} is equivalent to  $P=-D^{-2}R$ and there exists an $\eta_0$ such that $B\eta_0+DR\si-\rho=0$ and  $A\eta_0+q-D^{-2}R(b+C\si)=0$.
	
	If $\g\neq 0$,
	$\Psi[P]$ has a maximum $\b-2\sqrt{\a\g}$	
	and the maximum point is taken at  $P_*=\sqrt{\gamma/\a}-D^{-2}R>-D^{-2}R$.	 Then  Problem (EC) is finite if $\b-2\sqrt{\a\g}>0$.\ss
	
	If $\g=0$ and $\b-2\sqrt{\a\g}=0$, take $P_*=\sqrt{\gamma/\a}-D^{-2}R>-D^{-2}R$, the fifth line of \eqref{griccati-ass0000} is equivalent to
	$$q+P_*b+\Th_*\rho+(C+D\Th_*)P_*\si\in\sR(A+B\Th_*),$$
	where $$\ba{ll}\Theta_*\ad=-D^{-2}\sqrt{\frac{\a}{\g}}\[\(\sqrt{\frac{\gamma}{\a}}-D^{-2}R\)\(B+CD\)+S\]\\
	\ns\ad=-D^{-2}(B+CD)-D^{-2}\sqrt{\frac{\a}{\g}}\(S-D^{-2}R(B+CD)\)\\
	\ns\ad=-D^{-2}(B+CD)-|D^{-1}|\sqrt{\a}\cdot\text{sgn}(S-D^{-2}R(B+CD)).\ea$$
	We also notice that $|D^2\Th_*+(B+CD)|=|D|\sqrt{\a}$.
	Such $\Theta_*\notin \BTh[A,C;B,D].$\ss %(\textcolor{blue}{This case looks like all the quadratic term are canceled...thus the optimal one takes on the boundary of $\BTh$.})
	
(2) {\bf  Solvability}. Note that the sufficient condition (H3) for solvability requires ${\rm (LQ)}^0_\infty$	 to be solvable, we can directly take Theorem 7.2 (especially the (7.15) in the proof)  from \cite{Sun-Yong 2018} which presents an equivalence characterization. Then the result follows directly by verifying \eqref{H3equi-2}. If $\b-2\sqrt{\a\g}>0$, \eqref{H3equi-2} holds naturally.  If $\b=\g=0$, we have to take $P=-D^{-2}R$. Then \eqref{H3equi-2} is equivalent to  there exists an $\eta_0$ such that $ B\eta_0-D^{-1}R\si+\rho=A\eta_0+q-D^{-2}R(b+C\si)=0$.  The proof is complete.

\end{document}